\documentclass[a4paper,10pt]{amsart}
\usepackage{amssymb}
\usepackage{amsmath}
\usepackage{a4wide}
\newtheorem{theorem}{Theorem}
\def\afrac#1#2{#1/#2}

\begin{document}

\title{Metric discrepancy results for geometric progressions with small ratios}
\author{K. Fukuyama}
\address{Department of Mathematics, Kobe University, Rokko Kobe 657-8501 Japan}
\email{fukuyama@math.kobe-u.ac.jp}
\author{S. Sakaguchi}
\address{Aioi Nissay Dowa Insrance}
\author{O. Shimabe}
\address{Hamada Electrical Industries}
\author{T. Toyoda}
\address{Department of Mathematics, Kobe University, Rokko Kobe 657-8501 Japan}
\author{M. Tscheckl}
\address{Department of Mathematics, University of Graz}
\thanks{The research is partially supported by JSPS KAKENHI 16K05204}

\maketitle

\section{Introduction}

A sequence $\{x_k\}$ of real numbers is said to be uniformly distributed mod 1 if
$$
\lim_{N\to\infty}
\frac1N {}^\# \{ k\le N \mid \langle x_k \rangle \in [\,a, b)\}
= b-a\quad\hbox{for all}\quad0\le a< b< 1,
$$
where $\langle x\rangle$ denotes the fractional part $x-[\,x\,]$ of $x$. 
Since the convergence is uniform in $a$ and $b$, 
we use the following discrepancy $D_N(\{x_k\})$ to 
measure the speed of convergence: 
$$
D_N(\{x_k\})
= \sup _{0\le a < b\le 1} 
\biggl|
\frac1N {}^\# \{ k\le N \mid \langle x_k \rangle \in [\,a, b)\}
- (b-a )\biggr|. 
$$

For geometric progressions $\{\theta^k x\}$ with  $|\theta|>1$, 
we can prove the law of the iterated logarithm in exact form as below 
and determine the speed of convergence toward the uniform distribution. 
$$
\varlimsup _{N\to \infty}
\frac{ND_N (\{\theta^k x\})}{\sqrt {2N\log \log N}}
=
\Sigma_\theta
\quad \hbox{a.e. } x,
$$
where $\Sigma_\theta\ge 1/2$ is a constant determined by $\theta$. 
The case $\theta>1$ was proved in 
\cite{2008F} and the case $\theta< -1$ in \cite{geom}. 

Before this result, 
Philipp \cite{1975Ph} applied the method of Takahashi \cite{1962T} 
and proved that the limsup above is bounded from below and above by 
positive constants if we replace $\theta^k$ by $n_k$ satisfying $n_{k+1}/n_k \ge q> 1$. 

We are interested in the concrete value of $\Sigma_\theta$ because it indicates the 
speed of convergence. 

When $\theta^k \notin \bf Q$ for all $k=1$, $2$, \dots, then 
\begin{equation}\label{Eq:nonroot}
\Sigma_\theta=\frac12.
\end{equation}

When $\theta^k \in \bf Q$ for some $k=1$, $2$, \dots, 
 denote  $r=\min\{k\in {\bf N}: \theta^k \in \bf Q\}$ and
$\theta^r = p/q$ by $p\in \mathbf Z$ and $q\in \mathbf N$ with $\gcd(p,q)=1$. 
We first see that $\Sigma_\theta$ is independent of $r$ and 
and is determined only by $p$ and $q$, i.e., 
\begin{equation}\label{Eq:rindep}
\Sigma_\theta=\Sigma_{p/q}.
\end{equation}
We also have the estimate
\begin{equation}\label{Eq:sbound}
\frac12< \Sigma_{p/q}\le \Sigma_{|p|/q}\le 
\frac12\sqrt{\frac{|p|q+1}{|p|q-1}}.
\end{equation}
If both $p$ and $q$ are odd, then
\begin{equation}\label{Eq:type0}
\Sigma_{p/q} = \frac12\sqrt{\frac{|p|q+1}{|p|q-1}}.
\end{equation}
If  $|p|\ge 4$ is even and $q=1$, then 
\begin{equation}\label{Eq:even}
\Sigma_p=\frac12\sqrt{ \frac{(|p|+1)|p|(|p|-2)}{(|p|-1)^3}}
\end{equation}
If $p=2$ and $q=1$, then 
\begin{equation}\label{Eq:two}
\Sigma_2= \frac19{\sqrt{42}}
.
\end{equation}
By (\ref{Eq:nonroot}), (\ref{Eq:rindep}), (\ref{Eq:sbound}), and (\ref{Eq:two}), 
we see that $\Sigma_2$ is largest among $\Sigma_\theta$ ($|\theta|>1$) and 
that $\{2^k x\}$ is furthest from the uniform distribution a.e. 
These results are  proved in \cite{2008F,2010Fphv,geom}. 
Determining the concrete value of $\Sigma_{-2}$ is rather hard work \cite{mtwo}: 
\begin{equation}\label{Eq:mtwo}
\Sigma_{-2}= \frac1{49}{\sqrt{910}}
.
\end{equation}

When $q> 1$ and $pq$ is even, evaluation of the concrete value of $\Sigma_{p/q}$ needs 
very delicate estimate. We succeeded in giving 
the  closed formula below to have the concrete evaluation when $p/q$ is large. 
If $p$ is odd, $q$ is even and $|p|/q\ge 9/4$, 
or if $p$ is even, $q$ is odd and $|p|/q\ge 4$, then \cite{ymst}
\begin{equation}\label{Eq:large}
\Sigma_{p/q}=
\sqrt{
\frac{(|p|q)^{I}+1}{(|p|q)^{I}-1}
v\Bigl(\frac{|p|-q-1}{2(|p|-q)}\Bigr)
+
\frac{2(|p|q)^{I}}{(|p|q)^{I}-1}
\sum_{m=1}^{I-1}
\frac1{(|p|q)^{m}}
v\Bigl(q^m\frac{|p|-q-1}{2(|p|-q)}\Bigr)
},
\end{equation}
where 
$I=\min\{ n\in {\bf N}\mid q^n = \pm 1 \hbox{ \rm mod } |p|-q\}$ and
$v(x)=\langle x \rangle(1-\langle x\rangle)$.

Having these results, 
it is very natural to have a question if the formula (\ref{Eq:large}) is valid 
when $p/q$ is small. 
We already have  counterexamples $\Sigma_2$ and $\Sigma_{-2}$, 
since $I=1$ and the right hand side of the formula equals to 0 in these cases
which is different from the actual values  (\ref{Eq:two}) and (\ref{Eq:mtwo}). 
In this note, we make a conjecture and give a few affirmative examples for that. 
Because the case when $p/q$ is negative is very delicate and hard to be investigated, 
we restrict ourselves to the case of positive $p/q$. 

If $q$ is positive,  we \cite{2008F} have proved
\begin{align*}
\Sigma_{\theta}^2
&=
\sup_{0\le a \le 1}
\sigma_{\theta}^2(a)
=
\sup_{0\le a \le 1}
\biggl(
V(a ,a)
+
2\sum_{k=1}^\infty
\frac{
V(\langle p^k a \rangle,   \langle q^k a \rangle)}
{p^kq^k}
\biggr),
\end{align*}
where $V(x,\xi)= x\wedge \xi - x\xi$. 

When $pq$ is even and $p/q$ is large, the formula (\ref{Eq:large}) derived from
\begin{equation}\label{Eq:typei}
\Sigma_{p/q} = \sigma_{p/q}\Bigl(\frac{p-q-1}{2(p - q)}\Bigr).
\end{equation}

This equation holds since 
the point $x=(p-q-1)/2(p-q)$ is the maximal point of the function 
$V(x,x) + 2V(\langle px\rangle, \langle qx\rangle)/ pq$ 
and the remainder terms can be negligible when $p/q$ is large. 

When $p/q$ is small, then the remainder terms grow and 
validity of the formula (\ref{Eq:typei}) cannot be expected. 
Actually we can find several $p/q$ such that
\begin{equation}\label{Eq:typek}
\Sigma_{p/q} = \sigma_{p/q}\Bigl(\frac{n}{p^k - q^k}\Bigr) 
\end{equation}
holds for some $n\in {\mathbf N}$, 
and we say that $p/q$ is of type $k$. 
The formula (\ref{Eq:typei}) shows that $p/q$ is of type I if $pq$ is even and $p/q$ is large.

Since we have
$$\Sigma_{2} = \sigma_{2}\Bigl(\frac{1}{2^2 - 1^2}\Bigr), $$
we see that $2=2/1$ is of type II. 

We conjecture that there exists $k$ and $n$ such that (\ref{Eq:typek}) holds if $pq$ is even. 
We are now in a position to state our result. 

\begin{theorem}
Ratio $13/6$ is under the threshold of  validity of (\ref{Eq:large}), it is of 
Type I  and the concrete evaluation given by (\ref{Eq:large}) is as follows: 
$$
\Sigma_{13/6}=\sigma_{13/6}\Bigl(\frac{3}{13-6}\Bigr)=
\frac27\sqrt{\frac{237}{77}}.
$$
Ratios $4/3$, $8/3$, $10/3$, $12/5$, $17/8$ are of Type II and the concrete evaluations
are as follows:
\begin{align*}
\Sigma_{4/3}&=\sigma_{4/3}\Bigl(\frac{3}{4^2-3^2}\Bigr)=
\frac{18}{7}\sqrt{\frac{117609}{2985983}},
\\
\Sigma_{8/3}&=\sigma_{8/3}\Bigl(\frac{24}{8^2-3^2}\Bigr)=
\frac{2}{275}\sqrt{\frac{157667789263012683051319944222}{32159909742724829389686571}},
\quad
\\
\Sigma_{10/3}&=\sigma_{10/3}\Bigl(\frac{40}{10^2-3^2}\Bigr)=
\frac{6}{637}\sqrt{\frac{43479927170}{14877551}},
\\
\Sigma_{12/5}&=\sigma_{12/5}\Bigl(\frac{55}{12^2-5^2}\Bigr)
%\\&
=\frac{2}{119}
\sqrt{
\frac{\openup-2\jot
\begin{aligned}
&201428358526705877909161869599911562733573899\\&35684732999641319454307059909247205938048130
\end{aligned}
}
{\openup-2\jot
\begin{aligned}
&2245225770735455724008721112379267481599999\\&9999999999999999999999999999999999999999999
\end{aligned}
}
,}
\\
\Sigma_{17/8}&=\sigma_{17/8}\Bigl(\frac{101}{17^2-8^2}\Bigr)
%\\&
=
\frac{2}{675}\sqrt{ \frac{2972498069993207680191880231521354312204246283}
{104127550853833809985880737775289062764635}}.
\end{align*}
Ratio $19/10$ is of type III and the concrete evaluation is as follows:
\begin{align*}
\Sigma_{19/10}&=\sigma_{19/10}\Bigl(\frac{2879}{19^3-10^3}\Bigr)
%\\&
=
\frac{2}{17577}
\sqrt{\frac{
\openup-2\jot\begin{aligned}
&16497674472384448745389239618125365389\\&8900368036274961322838030028584448835
\end{aligned}
}
{
\openup-2\jot\begin{aligned}
&8535800662859082038722574792812344\\&200037037037037037037037037037037
\end{aligned}
}
.}
\end{align*}
Ratio $12/7$ is of type IV and the concrete evaluation is as follows:
\begin{align*}
\Sigma_{12/7}&= \sigma_{12/7}\Bigl(\frac{8717}{12^4-7^4}\Bigr)
\\&
=\frac{1}{18335}\sqrt{
\frac{1288914789424650371352900618359881195696318380071236938}
{15230103878098355389592475654267327331681959935}
}.
\end{align*}
Ratio $8/5$ is of type V and the concrete evaluation is as follows:
\begin{align*}
\Sigma_{8/5}&= \sigma_{8/5}\Bigl(\frac{13690}{8^5-5^5}\Bigr)
%\\&
=
\frac{1}{326073}
\sqrt{
\frac{
\openup-2\jot\begin{aligned}
&2693647024766931825274236270928683791\\&388436386344146949339630859940359610
\end{aligned}
}
{
\openup-2\jot\begin{aligned}
&9991122476153960121538644628099\\&1735537190082644628099173553719
\end{aligned}
}
}.
\end{align*}
Ratio $3/2$ is of type VI and the concrete evaluation is as follows:
\begin{align*}
\Sigma_{3/2}&= \sigma_{3/2}\Bigl(\frac{277}{3^6-2^6}\Bigr)
=
\frac{2}{665}\sqrt{\frac{305671451762616889661445636790873}{10314424798490535546171949055}}
.
\end{align*}
\end{theorem}
Since the evaluation of $\Sigma_{3/2}$ contains very delicate calculation 
and the proof is lengthy, 
it  will be proved in a separate paper \cite{aspm}. 
By the above results, 
concrete values of $\Sigma_{p/2}$ and $\Sigma_{p/3}$ are completely determined.

\section{Preliminary}
We denote $\sigma^2_\theta(a)$ simply by $\sigma^2(a)$. 

The inequality $0\le V(x,y) \le 1/4$ ($x$, $y\in [\,0,1)$) implies
\begin{align}\label{Eq:est}
\biggl|
2\sum_{n=N+1}^\infty
\frac{1}{(pq)^n}
V\bigl( \langle p^n x\rangle, \langle q^n x\rangle\bigr)
\biggr|
&\le 
\frac12\sum_{n=N+1}^\infty
\frac{1}{(pq)^n}
=
\frac{1}{2(pq-1)(pq)^N }
.
\end{align}

If $x> y$, then we have $V(x,y) = y(1-x)$, 
$V_x(x,y) = -y<0$, and $V_y(x,y) = 1-x> 0$. 
If $x< y$, then we have $V(x,y) = x(1-y)$, 
$V_x(x,y) = 1-y>0$, and $V_y(x,y) = -x< 0$. 
Since one of  $V_x(x,y)$ and  $V_y(x,y)$ is positive and the other is negative, 
we see
$$
\Bigl|
\frac{d}{dx} 
V\bigl( \langle p^n x\rangle, \langle q^n x\rangle\bigr)
\Bigr|
=
\bigl|
p^n V_x\bigl( \langle p^n x\rangle, \langle q^n x\rangle\bigr)
+
q^n V_y\bigl( \langle p^n x\rangle, \langle q^n x\rangle\bigr)
\bigr|
\le p^n
\quad\hbox{a.e.}
$$
and 
\begin{align}\label{Eq:dest}
\biggl|
2\sum_{n=N+1}^\infty
\frac{1}{(pq)^n}
\frac{d}{dx} 
V\bigl( \langle p^n x\rangle, \langle q^n x\rangle\bigr)
\biggr|
&\le
2\sum_{n=N+1}^\infty
\frac{1}{q^n}
=
\frac{2}{(q-1)q^N}
\qquad\hbox{a.e.}
\end{align}
Since the series (\ref{Eq:est}) and (\ref{Eq:dest}) are uniformly convergent, 
we see 
$$
\sigma^2(a)
=
\int_0^a 
\biggl(
\frac{d}{dx}V(x ,x)
+
2\sum_{k=1}^\infty
\frac{1}{p^kq^k}
\frac{d}{dx}
V(\langle p^k a \rangle,   \langle q^k a \rangle)
\biggr)
\,dx.
$$
It make us possible to conclude that $\sigma^2(a)$ increases in the interval 
in which 
$$\frac{d}{dx}V(x ,x)
+
2\sum_{k=1}^\infty
\frac{1}{p^kq^k}
\frac{d}{dx}
V(\langle p^k a \rangle,   \langle q^k a \rangle)
\ge 0
\quad\hbox{a.e.}
$$

\section{Type IV case, $\Sigma_{12/7}$}

Put $
c = \dfrac{8717}{12^4-7^4}= \dfrac{8717}{18335}$. 
By 
$$
12^{24}= 7^{24}=1 \quad\hbox{mod}\quad 12^4-7^4
,$$
we have $V\bigl(\langle 7^{n+24k}c\rangle, \langle 12^{n+24k}c\rangle\bigr)
=V\bigl(\langle 7^{n}c\rangle, \langle 12^{n}c\rangle\bigr)$
and
\begin{align*}
\sigma^2(c)
&=
V\bigl(\langle c\rangle, \langle c\rangle\bigr)
+
2
\frac{7^{24}12^{24}}{7^{24}12^{24}-1}
\sum_{n=1}^{24}
\frac1{7^n12^n}
V\bigl(\langle 7^{n}c\rangle, \langle 12^{n}c\rangle\bigr)
\\&=
\frac{1288914789424650371352900618359881195696318380071236938}
{5119937907681452900160044383953378173894837463709805375}
.
\end{align*}
We divide $[\,0,1/2)$ into 
 $[\,0,3/7)$,
 $[\,3/7,68/12^2)$, 
 $[\,68/12^2, 45/95)$,
 $[\,45/95, 820/12^3)$, 
 $[\,820/12^3,821/12^3)$,
 $[\,821/12^3,9858/12^4)$,  
 $[\,9858/12^4, c)$, 
 $[\,c, 9859/12^4)$, 
 $[\,9859/12^4, 822/12^3)$, 
 $[\,822/12^3,659/1385)$, \break
 $[\,659/1385, 69/12^2)$, 
 $[\,69/12^2,46/95)$, 
 $[\,46/95,1/2)$, 
and prove $\sigma^2(x) < \sigma^2(c) $ ($x\not=c$) on each.

\subsection{ $[\,0,3/7)$ part.}
On $[\,0,3/7)$, by applying (\ref{Eq:est}) for $N=0$   we have
\begin{align*}
\sigma^2(x) -\sigma^2(c) 
&
\le 
x(1-x) + \frac{1}{2\cdot 83} - \sigma^2(c)\biggl|_{x=3/7}
\\&
=
-\frac{412525633762177727202928217956973408126390844189499799}
{501753914952782384215684349627431061041694071443560926750}
<0.
\end{align*}

\subsection{ $[\,3/7,68/12^2)$ and   $[\,46/95,1/2)$ parts.}
On $[\,3/7, 1/2)$, we see
$\langle 12x\rangle = 12x-5$, 
$\langle 7x\rangle = 7x-3$, 
$\langle 12x\rangle -\langle 7x\rangle = 5x-2> 0$, 
and 
\begin{equation}\label{Eq:12/7second}
V\bigl( \langle 7x\rangle , \langle 12x\rangle \bigr)=
(7x-3) (6-12x)
\quad
x\in [\,3/7, 1/2)
.
\end{equation}

On $[\,3/7,68/12^2)$, by applying (\ref{Eq:est}) for $N=1$  we have
\begin{align*}
\sigma^2(x) -\sigma^2(c)
&
\le
x(1-x) + \frac{2}{7\cdot 12}(7x-3)(6-12x) + \frac{1}{2\cdot 83\cdot 7\cdot 12}-\sigma^2(c)
\biggl|_{x=68/12^2}
\\
&=
-\frac{295227494577955681873105793071261028441169351695037887}
{15482692232828713570083974217075015597857988490258451454000}
< 0, 
\end{align*}
since the bounding quadratic function on the right hand side 
has the axis of symmetry at $x=10/21\in (68/12^2, 46/95)$. 
On $[\,46/95,1/2)$, it also implies
\begin{align*}
\sigma^2(x) -\sigma^2(c)
&
\le
x(1-x) + \frac{2}{7\cdot 12}(7x-3)(6-12x) + \frac{1}{2\cdot 83\cdot 7\cdot 12}-\sigma^2(c)
\biggl|_{x=46/95}
\\
&=
-\frac{141744667266329259543498161264701235944975820204503619}
{860149568490484087226887456504167533214332693903247303000}
< 0.
\end{align*}

\subsection{ $[\,69/12^2,46/95)$ part.}
We have  $\langle 12^2 x\rangle = 12^2x -69$, 
$\langle 7^2 x\rangle = 7^2 x -23$, and
$\langle 12^2 x\rangle -\langle 7^2 x\rangle = 95x-46< 0$. 
By applying (\ref{Eq:est}) for $N=2$   we have
\begin{align*}
\sigma^2(x) -\sigma^2(c)
&
\le
x(1-x) + \frac{2}{7\cdot 12}(7x-3)(6-12x) + 
\frac{2}{7^2\cdot 12^2}(12^2x-69)(24-7^2x)
\\&\qquad\qquad
+
\frac{1}{2\cdot 83\cdot 7^2\cdot 12^2}-\sigma^2(c)
\biggl|_{x=5639/11760}
\\
&=
-\frac{648440595936525455567931682013588012496682952191224629959}
{9441001663751553341402316722589742844560515648282042397728000}
< 0
,
\end{align*}
since the quadratic function has the axis at $x=5639/11760$.

\subsection{ $[\,68/12^2, 45/95)$, $[\,45/95, 820/12^3)$, 
and  $[\,659/1385, 69/12^2)$ parts.}

On $[\,68/12^2, 69/12^2)$, 
we see $\langle 12^2 x\rangle = 12^2x -68$, 
$\langle 7^2 x\rangle = 7^2 x -23$, 
$\langle 12^2 x\rangle -\langle 7^2 x\rangle = 95x-45$
and 
\begin{equation}\label{Eq:12/7third}
V\bigl(\langle 7^2 x\rangle, \langle 12^2 x\rangle \bigr)
=
\begin{cases}
(12^2x -68)(24-7^2x)
& x\in [\,68/12^2, 45/95),
\\
(7^2x-23)(69-12^2x)
& x\in [\,45/95,69/12^2).
\end{cases}
\end{equation}

On $[\,68/12^2, 45/95)$, by applying (\ref{Eq:est}) for $N=2$  we have
\begin{align*}
\sigma^2(x) -\sigma^2(c)
&
\le
x(1-x) + \frac{2}{7\cdot 12}(7x-3)(6-12x) + 
\frac{2}{7^2\cdot 12^2}(12^2x-68)(24-7^2x)
\\&\qquad\qquad
+
\frac{1}{2\cdot 83\cdot 7^2\cdot 12^2}-\sigma^2(c)
\biggl|_{x=45/95}
\\
&=
-\frac{1042072175578582511796226709992117624556912690252377931}
{72252563753200663327058546346350072790003946287872773452000}
<0
,\end{align*}
since the quadratic function  has the axis at
$x={4217}/{8820}> 45/95$. 

On $[\,45/95, 820/12^3)$, by applying (\ref{Eq:est}) for $N=2$  we have
\begin{align*}
\sigma^2(x) -\sigma^2(c)
&
\le
x(1-x) + \frac{2}{7\cdot 12}(7x-3)(6-12x) + 
\frac{2}{7^2\cdot 12^2}(7^2x-23)(69-12^2x)
\\&\qquad\qquad+
\frac{1}{2\cdot 83\cdot 7^2\cdot 12^2}-\sigma^2(c)
\biggl|_{x=820/12^3}
\\
&=
-\frac{150376946989927311461065482574253514791344687287872790913}
{46819661312074029835933938032434847167922557194541557196896000}
<0, \end{align*}
since the quadratic function  has the axis at
$x={5591}/{11760}\in (820/12^3, 659/1385)$. 
On $[\,659/1385,69/12^2)$ it also implies
\begin{align*}
\sigma^2(x) -\sigma^2(c)
&
\le
x(1-x) + \frac{2}{7\cdot 12}(7x-3)(6-12x) + 
\frac{2}{7^2\cdot 12^2}(7^2x-23)(69-12^2x)
\\&\qquad\qquad+
\frac{1}{2\cdot 83\cdot 7^2\cdot 12^2}-\sigma^2(c)
\biggl|_{x=659/1385}
\\
&=
-\frac{95246882859225298205690277732878224974396743417692962899}
{5543866964219333696421875202609094735104212794722190034198508000}
< 0.\end{align*}

\subsection{ $[\,820/12^3,821/12^3)$, $[\,821/12^3,9858/12^4)$, 
 $[9859/12^4, 822/12^3)$, and $[\,822/12^3,659/1385)$ parts.}

On $[\,820/12^3, 821/12^3)$ we have
$\langle 12^3 x\rangle = 12^3 x-820$, 
$\langle 7^3 x\rangle = 7^3x-162$, 
$\langle 12^3 x\rangle -\langle 7^3 x\rangle = 1385 x-658$, and 
$$
V\bigl( \langle 7^3 x\rangle, \langle 12^3 x\rangle\bigr)
=
\begin{cases}
(12^3x-820)(163-7^3x) & x\in [\,820/12^3, 658/1385),
\\
(7^3x-162)(821-12^3x) & x\in [\,658/1385, 821/12^3).
\end{cases}
$$

On $[\,820/12^3, 658/1385)$, by applying (\ref{Eq:est}) for $N=3$  
we have
\begin{align*}
\sigma^2(x) -\sigma^2(c)
&
\le
x(1-x) + \frac{2}{7\cdot 12}(7x-3)(6-12x) + 
\frac{2}{7^2\cdot 12^2}(7^2x-23)(69-12^2x)
\\&\qquad+
\frac{2}{7^3\cdot 12^3}(12^3x-820)(163-7^3x)
+
\frac{1}{2\cdot 83\cdot 7^3\cdot 12^3}-\sigma^2(c)
\biggl|_{x=658/1385}
\\
&=
-\frac{82317377154826056601317125454132578280311183651631748395617}
{155228274998141343499812505673054652582917958252221320957558224000}
< 0, 
\end{align*}
since the quadratic function has the axis at 
$x=123241/259308> 658/1385$. 

On $[\,658/1385,821/12^3)$,  by applying (\ref{Eq:est}) for $N=3$  we have
\begin{align*}
\sigma^2(x) -\sigma^2(c)
&
\le
x(1-x) + \frac{2}{7\cdot 12}(7x-3)(6-12x) + 
\frac{2}{7^2\cdot 12^2}(7^2x-23)(69-12^2x)
\\&\qquad+
\frac{2}{7^3\cdot 12^3}(7^3x-162)(821-12^3x)
+
\frac{1}{2\cdot 83\cdot 7^3\cdot 12^3}-\sigma^2(c)
\biggl|_{x=658/1385}
\\
&=
-\frac{82317377154826056601317125454132578280311183651631748395617}
{155228274998141343499812505673054652582917958252221320957558224000}
< 0, 
\end{align*}
since the quadratic function 
has the axis at 
$x=1970471/4148928< 658/1385$.

On $[\, 821/12^3, 163/7^3)$, 
we see $\langle 12^3 x\rangle = 12^3 x-821$, 
$\langle 7^3 x\rangle = 7^3x-162$, and
$\langle 12^3 x\rangle -\langle 7^3 x\rangle = 1385 x-659< 0$.
By applying (\ref{Eq:est}) for $N=3$   we have
\begin{align*}
\sigma^2(x) -\sigma^2(c)
&
\le
x(1-x) + \frac{2}{7\cdot 12}(7x-3)(6-12x) + 
\frac{2}{7^2\cdot 12^2}(7^2x-23)(69-12^2x)
\\&\qquad+
\frac{2}{7^3\cdot 12^3}(12^3x-821)(163-7^3x)
+
\frac{1}{2\cdot 83\cdot 7^3\cdot 12^3}-\sigma^2(c)
\biggl|_{x=163/7^3}
\\
&=
-\frac{676061216533256609810614263080777776833310040575587371597}
{2081740866857217511779210837331038297225593700446190348699024000}
< 0, 
\end{align*}
since the quadratic function has the axis at 
$x=1972199/4148928>163/7^3$. 

On $[\,  163/7^3, 822/12^3)$, 
we see $\langle 12^3 x\rangle = 12^3 x-821$, 
$\langle 7^3 x\rangle = 7^3x-163$, 
$\langle 12^3 x\rangle -\langle 7^3 x\rangle = 1385 x-658> 0$ 
by $658/1385< 163/7^3$, 
and
\begin{equation}\label{Eq:12/7fourth}
V\bigl( \langle 7^3 x\rangle, \langle 12^3 x\rangle\bigr)
=
(7^3x-163)(822-12^3x)
\quad
x\in [\,  163/7^3, 822/12^3).
\end{equation}

On $[\,163/7^3, 9858/12^4)$, by applying (\ref{Eq:est}) for $N=3$  we have
\begin{align*}
\sigma^2(x) -\sigma^2(c)
&
\le
x(1-x) + \frac{2}{7\cdot 12}(7x-3)(6-12x) + 
\frac{2}{7^2\cdot 12^2}(7^2x-23)(69-12^2x)
\\&\qquad+
\frac{2}{7^3\cdot 12^3}(7^3x-163)(822-12^3x)
+
\frac{1}{2\cdot 83\cdot 7^3\cdot 12^3}-\sigma^2(c)
\biggl|_{x=9858/12^4}
\\
&=
-\frac{118487834091551690955153210037952312413729819138716728399}
{20975208267809165366498404238530811531229305623154617624209408000}
< 0, 
\end{align*}
since the quadratic function has the axis at 
$x=328757/691488\in [\,9858/12^4,9859/12^4)$. 
On $[\,9859/12^4,822/12^3)$, it also implies
\begin{align*}
\sigma^2(x) -\sigma^2(c)
&
\le
x(1-x) + \frac{2}{7\cdot 12}(7x-3)(6-12x) + 
\frac{2}{7^2\cdot 12^2}(7^2x-23)(69-12^2x)
\\&\qquad+
\frac{2}{7^3\cdot 12^3}(7^3x-163)(822-12^3x)
+
\frac{1}{2\cdot 83\cdot 7^3\cdot 12^3}-\sigma^2(c)
\biggl|_{x=9859/12^4}
\\
&=
-\frac{1751672514624267500406933768845174563511908294312287783239}
{755107497641129953193942552587109215124255002433566234471538688000}
<0. 
\end{align*}

On $[\,822/12^3,659/1385)$, 
we see $\langle 12^3 x\rangle = 12^3 x-822$, 
$\langle 7^3 x\rangle = 7^3x-163$, and
$\langle 12^3 x\rangle -\langle 7^3 x\rangle = 1385 x-659< 0$. 
By applying (\ref{Eq:est}) for $N=3$  we have
\begin{align*}
\sigma^2(x) -\sigma^2(c)
&
\le
x(1-x) + \frac{2}{7\cdot 12}(7x-3)(6-12x) + 
\frac{2}{7^2\cdot 12^2}(7^2x-23)(69-12^2x)
\\&\qquad+
\frac{2}{7^3\cdot 12^3}(12^3x-822)(164-7^3x)
+
\frac{1}{2\cdot 83\cdot 7^3\cdot 12^3}-\sigma^2(c)
\biggl|_{x=659/1385}
\\
&=
-\frac{146042660849435492820349042393734530006259332002143045874331}
{465684824994424030499437517019163957748753874756663962872674672000}
< 0, 
\end{align*}
since the quadratic function  has the axis at 
$x=329045/691488> 659/1385$.

\subsection{ $[\,9858/12^4, c)$ and $[\,c, 9859/12^4)$ parts.}
We see
$\langle 12^4 x\rangle=12^4x-9858$, 
$\langle 7^4 x\rangle=7^4x-1141$, 
and
$$
V\bigl(\langle 7^4 x\rangle, \langle 12^4 x\rangle\bigr)
=
\begin{cases}
(12^4x-9858)(1142-7^4x) 
& x\in [\,9858/12^4, c),
\\
(7^4x-1141)(9859-12^4x)
& x\in [\,c, 9859/12^4).
\end{cases}
$$
On $[\,9858/12^4, c)$, 
by applying (\ref{Eq:dest}) for $N=4$  we have 
\begin{align*}
\frac{d}{dx}
\sigma^2(x)
&\ge
\frac{d}{dx}
\Bigl(
x(1-x)
+
\frac{2}{7\cdot 12}(7x-3)(6-12x)
+
\frac{2}{7^2 12^2}(7^2x-23)(69-12^2x)
\\&\qquad\qquad
+
\frac{2}{7^3 12^3}(7^3x-163)(822-12^3x)
+
\frac{2}{7^4 12^4}(12^4x-9858)(1142-7^4x) 
\Bigr)
-\frac{2}{6\cdot 7^4}
\\
&
=
-\frac{74680704 x - 35506607}{4148928}
\ge
-\frac{74680704 c - 35506607}{4148928}
=\frac{21942577}{76070594880}
> 0.
\end{align*}
On $[\,c, 9859/12^4)$, by applying (\ref{Eq:dest}) for $N=4$   we have
\begin{align*}
\frac{d}{dx}
\sigma^2(x)
&\le
\frac{d}{dx}
\Bigl(
x(1-x)
+
\frac{2}{7\cdot 12}(7x-3)(6-12x)
+
\frac{2}{7^2 12^2}(7^2x-23)(69-12^2x)
\\&\qquad\qquad
+
\frac{2}{7^3 12^3}(7^3x-163)(822-12^3x)
+
\frac{2}{7^4 12^4}(7^4x-1141)(9859-12^4x) 
\Bigr)
+\frac{2}{6\cdot 7^4}
\\
&
=
-\frac{448084224 x - 213028219}{24893568}
\le 
-\frac{448084224 c - 213028219}{24893568}
=-\frac{77785243}{456423569280}
< 0.
\end{align*}

\section{Type I case, $\Sigma_{13/6}$}
Put
$$c = \frac37.$$
Because of $V(\langle 6^n c\rangle , \langle 13^n c\rangle) = V(c , c)$, 
we have
$$
\sigma^2(c) = \frac{6\cdot 13+1}{6\cdot 13-1}V(c,c) = 
\frac{948}{3773}
.$$
We divide $[\,0,1/2)$ into 
$[\,0,\afrac{5}{13})$, 
$[\,\afrac{5}{13}, c)$, 
$[\,c, \afrac{942}{13^3})$, 
$[\,\afrac{942}{13^3}, \afrac{943}{13^3})$,
$[\,\afrac{943}{13^3}, \afrac{73}{13^2})$, 
$[\,\afrac{73}{13^2}, \afrac{58}{133})$, 
$[\,\afrac{58}{133}, \afrac{6}{13})$, and
$[\,\afrac{6}{13}, \afrac12)$, 
and prove $\sigma^2(x) < \sigma^2(c) $ ($x\not=c$) on each. 

\subsection{ $[\,0,\afrac{5}{13})$ part.}
By applying (\ref{Eq:est}) for $N=0$, we have
\begin{align*}
\sigma^2(x) - \sigma^2(c) 
& \le 
\Bigl(x(1-x) + \frac{1}{2\cdot 77} - \sigma^2(c) \Bigr)_{x=5/13}
=
-\frac{10303}{1275274}
< 0
,
\end{align*}
since the quadratic function above has the axis at $x=1/2$.

\subsection{ $[\,\afrac{6}{13}, \afrac12)$ part.}
We have
$\langle 13x\rangle = 13x -6$, 
$\langle 6x\rangle = 6x -2$, and
$\langle 13x\rangle - \langle 6x\rangle = 7x -4<0$. 
By applying (\ref{Eq:est}) for $N=1$ we have
\begin{align*}
\sigma^2(x) - \sigma^2(c) 
& \le 
\Bigl(x(1-x) + \frac{2}{6\cdot 13}(13x-6)(3-6x)
+\frac{1}{2\cdot 77\cdot 6\cdot 13} - \sigma^2(c) \Bigr)_{x=19/39}
\\
& =
-\frac{1741}{2550548}
< 0
,
\end{align*}
since the quadratic function above has the axis at $x=\afrac{19}{39}$.

\subsection{ $[\,\afrac{58}{133}, \afrac{6}{13})$ part.}
On $[\,5/13, 6/13)$, we have
$\langle 13x \rangle = 13x-5$, 
$\langle 6x \rangle = 6x-2$, 
$\langle 13x \rangle -\langle 6x \rangle = 7x-3$, 
$3/7\in [\,5/13, 6/13)$, and
\begin{equation}\label{Eq:13/6second}
V\bigl(\langle 6x \rangle , \langle 13x \rangle \bigr)
=
\begin{cases}
(13x-5)(3-6x) & x\in [\,\afrac{58}{133}, c), 
\\
(6x-2)(6-13x) & x\in [\,c, \afrac{6}{13}). 
\end{cases}
\end{equation}
On $[\,\afrac{58}{133}, \afrac{6}{13})\subset[\,5/13, 6/13)$, 
by applying (\ref{Eq:est}) for $N=1$ we have
\begin{align*}
\sigma^2(x) - \sigma^2(c)
&\le 
\Bigl(x(1-x) + \frac2{6\cdot 13}(6x-2)(6-13x)
+
\frac1{2\cdot 6\cdot 13(6\cdot 13-1)}
-
\frac{948}{3773}\Bigr)\Bigm|_{x=58/133}
\\&
=
-\frac{6415}{212480268}< 0
,
\end{align*}
since the quadratic function has the axis at 
$x=\afrac{101}{234}< \afrac{58}{133}$.

\subsection{ $[\,\afrac{73}{13^2}, \afrac{58}{133})$ part.}
We have
$\langle 13^2 x\rangle = 13^2 x -73$, 
$\langle 6^2 x\rangle = 6^2 x -15$, and
$\langle 13^2 x\rangle -\langle 6^2 x\rangle  = 133x-58< 0$.
By applying (\ref{Eq:est}) for $N=2$ we have
\begin{align*}
&
\sigma^2(x) - \sigma^2(c)
\\&
\le 
\Bigl(x(1-x) + \frac{2}{6\cdot 13} (6x-2)(6-13x)
+ \frac{2}{6^2\cdot 13^2} (13^2 x-73)(16-6^2 x)
+
\frac{1}{2\cdot 6^2 \cdot 13^2 (6\cdot 13-1)}
\\&\qquad\qquad\qquad
- \sigma^2(c)\Bigr)\Bigm|_{x=1321/3042}
\\&
=-\frac{1842013}{69828903144}< 0,
\end{align*}
since the  quadratic function has the axis at $x=1321/3042$. 

\subsection{ $[\,\afrac{943}{13^3}, \afrac{73}{13^2})$ part.}
On $[\,72/13^2, 73/13^2)$, we have
$\langle 6^2 x \rangle = 6^2 x -15$, 
$\langle 13^2 x \rangle = 13^2 x -72$, 
$\langle 13^2 x \rangle- \langle 6^2 x \rangle  = 133(x-c)$, and
\begin{equation}\label{Eq:13/8third}
V\bigl(\langle 6^2 x \rangle, \langle 13^2 x \rangle\bigr)
=
(6^2 x -15)(73-13^2x) 
\qquad x\in  [\,c, 73/13^2) .
\end{equation}
On $[\,\afrac{943}{13^3}, \afrac{73}{13^2})\subset[\,c, 73/13^2)$, 
by applying (\ref{Eq:est}) for $N=2$ we have
\begin{align*}
\sigma^2(x) - \sigma^2(c)
&
\le 
\biggl(x(1-x) + \frac2{6\cdot 13}(6x-2)(6-13x)
+ \frac2{6^2\cdot 13^2}(6^2 x-15)(73-13^2x)
\\&\qquad\qquad
+ \frac1{2\cdot 6^2\cdot 13^2(6\cdot 13-1)}
-\sigma^2(c)
\biggr)\biggm|_{x=943/13^3}
\\&
=
-\frac{1662359}{1311231625704}< 0
,
\end{align*}
since the quadratic function has the axis at $x= 1449/3380< 943/13^3$.

\subsection{ $[\,\afrac{942}{13^3}, \afrac{943}{13^3})\subset [\,c, 73/13^2)$ part.}
We have
$\langle 6^3 x \rangle = 6^3 x -92$, 
$\langle 13^3 x \rangle = 13^3 x -942$,
and $\langle 13^3 x \rangle - \langle 6^3 x \rangle  = 1981x-850$. 

On $[\,\afrac{942}{13^3}, \afrac{850}{1981})$ 
by applying (\ref{Eq:est}) for $N=3$ we have
\begin{align*}
&\sigma^2(x)- \sigma^2(c)
\\
&
\le
\Bigl(
x(1-x) + \frac{2}{6\cdot 13} (6x-2)(6-13x)
+ \frac{2}{6^2\cdot 13^2} (6^2 x-15)(73-13^2)
\\&\qquad
+ \frac{2}{6^3\cdot 13^3} (13^3 x -942)(93-6^3x)
+
\frac{1}{2\cdot 6^3 \cdot 13^3 (6\cdot13-1)}
-\sigma^2(c)
\Bigr)\Bigm|_{x=474997/1107288}
\\&
=
-\frac{63760513}{94408677050688}
,
\end{align*}
since the quadratic function has the axis at $x=474997/1107288$.

On $[\, \afrac{850}{1981}, \afrac{943}{13^3})$, 
by applying (\ref{Eq:est}) for $N=3$ we have
\begin{align*}
\sigma^2(x)- \sigma^2(c)
&
\le
\Bigl(
x(1-x) + \frac{2}{6\cdot 13} (6x-2)(6-13x)
+ \frac{2}{6^2\cdot 13^2} (6^2 x-15)(73-13^2)
\\&\qquad
+ \frac{2}{6^3\cdot 13^3} (6^3x-92)(943-13^3x)
+
\frac{1}{2\cdot 6^3 \cdot 13^3 (6\cdot13-1)}
-\sigma^2(c)
\Bigr)\Bigm|_{x=850/1981}
\\&
=
-\frac{16535003}{22061250587376}< 0
,
\end{align*}
since the quadratic function has the axis at $x=711505/1660932< 850/1981$.

\subsection{ $[\,\afrac{5}{13}, c)$  part.}
By applying (\ref{Eq:dest}) for $N=1$ we have
\begin{align*}
\frac{d}{dx}\sigma^2(x) 
&\ge
\frac{d}{dx}\Bigl(
x(1-x) +\frac{2}{6\cdot 13}(13x-5)(3-6x)\Bigr) - \frac{2}{6\cdot 5}
=\frac{527}{195}-6x
> \frac{527}{195}-6c
=\frac{179}{1365}
> 0.
\end{align*}

\subsection{ $[\,c, \afrac{942}{13^3})$ part.}
We have $\langle 6^3x\rangle = 6^3x-92< \langle 13^3x\rangle = 13^3 x -941$. 
By applying (\ref{Eq:dest}) for $N=3$ we have
\begin{align*}
\frac{d}{dx}
\sigma^2(x)
&
\le
\frac{d}{dx}
\Bigl(
x(1-x) + \frac{2}{6\cdot 13}(6 x-2)(6-13 x)
+ \frac{2}{6^2\cdot 13^2} (6^2 x - 15)(73-13^2 x)
\\&\qquad\qquad\qquad
+ \frac{2}{6^3\cdot 13^3} (6^3 x-92)(942-13^3 x)
\Bigr)
+
\frac1{540}
\\&
=
-\frac{16609320 x - 7116167}{1186380}
\le -\frac{16609320 c - 7116167}{1186380}
=
-\frac{2113}{1186380}
<0
.
\end{align*}

\section{Type II case, $\Sigma_{4/3}$}

Put
$
c = \dfrac{3}{7}
$. 
By $3^{6}=4^{6}=1$ mod 7, 
we see
$V\bigl(\langle 3^{k+6j}c\rangle, \langle 4^{k+6j}c\rangle\bigr) 
=V\bigl(\langle 3^{k}c\rangle, \langle 4^{k}c\rangle\bigr) 
$
and
\begin{align*}
\sigma^2(c)&=
V(c, c)
+
\frac{2(3\cdot 4)^{6}}{(3\cdot 4)^{6}-1}\sum_{k=1}^{6}\frac1{(3\cdot4)^{k}} 
V\bigl(\langle 3^{k}c\rangle, \langle 4^{k}c\rangle\bigr) 
%\\&
=
\frac{38105316}{146313167}
.
\end{align*}
We divide $[\,0,1/2)$ into 
$[\,0,{1}/{4})$, 
$[\,{1}/{4}, {1}/{3})$, 
$[\,{1}/{3}, {3}/{8})$, 
$[\,3/{8}, {27}/{64})$, 
$[\,{27}/{64}, c)$, 
$[\,c, {16}/{37})$, 
$[\,{16}/{37}, {7}/{16})$
$[\,{7}/{16}, {4}/{9})$,
$[\,{4}/{9}, {1}/{2})$, 
and prove $\sigma^2(x) <  \sigma^2(c) $ $(x\not=c)$ on each. 

\subsection{ $[\,0,{1}/{4})$ part.}
On $[\,0,{1}/{4})$, by applying (\ref{Eq:est}) for $N=0$ we have
\begin{align*}
\sigma^2(x) - \sigma^2(c) 
&\le 
x(1-x) \Bigr|_{x={1}/{4}}
+ \frac1{22}- \sigma^2(c) 
=
-\frac{64335979}{2341010672}
<0.
\end{align*}

\subsection{ $[\,{1}/{4}, {1}/{3})$ part.}

We have
$\langle 4x\rangle = 4x -1$, 
$\langle 3x\rangle = 3x $, and
$\langle 4x\rangle - \langle 3x\rangle = x-1<0$.
By applying (\ref{Eq:est}) for $N=1$ we have
\begin{align*}
\sigma^2(x) - \sigma^2(c) 
&\le 
\Bigl(x(1-x) + \frac{2}{3\cdot 4} (4x-1)(1-3x)
\Bigr)\Bigr|_{x={1}/{3}}
 + \frac1{22\cdot 12}- \sigma^2(c) 
=
-\frac{362668489}{10534548024}
<0
,
\end{align*}
since the quadratic function has the axis  
at $x={13}/{36} > {1}/{3}$.

\subsection{ $[\,{1}/{3}, {3}/{8})$ part.}

On $[\,1/3, 1/2)$, we have
$\langle 4x\rangle = 4x -1$, 
$\langle 3x\rangle = 3x -1$, 
$\langle 4x\rangle - \langle 3x\rangle = x>0$, to have
\begin{equation}\label{Eq:4/3second}
V\bigl( \langle 4x\rangle , \langle 3x\rangle\bigr)
=
(3x-1)(2-4x)
\quad
x\in[\,1/3, 1/2)
.
\end{equation}

On $[\,1/{3}, {3}/{8})$, 
by applying (\ref{Eq:est}) for $N=1$ 
we have
\begin{align*}
\sigma^2(x) - \sigma^2(c) 
&\le 
\Bigl(x(1-x) + \frac{2}{3\cdot 4} (3x-1)(2-4x)
\Bigr)\Bigr|_{x={3}/{8}}
 + \frac1{22\cdot 12}- \sigma^2(c) 
\\&
=
-\frac{111030749}{9364042688}
<0
,
\end{align*}
since the quadratic function has the axis  
at $x={4}/{9} > {3}/{8}$.

\subsection{ $[\,{4}/{9}, {1}/{2})$ part.}

On 
$[\,4/{9}, {1}/{2})$, 
we have
$\langle 4^2 x\rangle = 4^2 x -7$, 
$\langle 3^2 x\rangle = 3^2 x -4$, and
$\langle 4^2 x\rangle - \langle 3^2 x\rangle = 7x -3 > 0$.
By applying (\ref{Eq:est}) for $N=2$ we have
\begin{align*}
&\sigma^2(x) - \sigma^2(c) 
\\
&\le 
\Bigl(x(1-x) + \frac{2}{3\cdot 4} (3x-1)(2-4x)
+ \frac{2}{3^2\cdot 4^2} (3^2x-4)(8-4^2x)
\Bigr)\Bigr|_{x={41}/{90}}
 + \frac1{22\cdot 12^2}- \sigma^2(c) 
\\&
=
-\frac{463595039}{1896218644320}
<0
,
\end{align*}
since the quadratic function has the axis  
at $x={41}/{90} $.

\subsection{ $[\,{7}/{16}, {4}/{9})$ part.}

On 
$[\,7/{16}, {4}/{9})$, 
we have
$\langle 4^2 x\rangle = 4^2 x -7$, 
$\langle 3^2 x\rangle = 3^2 x -3$, and
$\langle 4^2 x\rangle - \langle 3^2 x\rangle = 7x -4 < 0$.
By applying (\ref{Eq:est}) for $N=2$ we have 
\begin{align*}
&\sigma^2(x) - \sigma^2(c) 
\\
&\le 
\Bigl(x(1-x) + \frac{2}{3\cdot 4} (3x-1)(2-4x)
+ \frac{2}{3^2\cdot 4^2} (4^2 x-7)(4-3^2 x) 
\Bigr)\Bigr|_{x={319}/{720}}
 + \frac1{22\cdot 12^2}- \sigma^2(c) 
\\&
=
-\frac{12926489833}{15169749154560}
<0
,
\end{align*}
since the quadratic function has the axis  
at $x={319}/{720} $.

\subsection{ $[\,3/{8}, {27}/{64})$ and  $[\,{16}/{37}, {7}/{16})$ parts.}

On 
$[\,3/{8}, {7}/{16})$, 
we have
$\langle 4^2 x\rangle = 4^2 x -6$, 
$\langle 3^2 x\rangle = 3^2 x -3$, 
$\langle 4^2 x\rangle - \langle 3^2 x\rangle = 7x -3 < 0$, 
and
\begin{equation}\label{Eq:4/3third}
V\bigl(\langle 3^2 x\rangle, \langle 4^2 x\rangle\bigr)
=
\begin{cases}
(4^2 x -6)(4-3^2 x) & x\in [\,3/{8}, c),\\
(3^2 x-3)(7-4^2 x) & x\in [\,c, {7}/{16}).
\end{cases}
\end{equation}

On $[\,3/{8}, {27}/{64})$, by noting $27/64< c$ and 
by applying (\ref{Eq:est}) for $N=2$ we have
\begin{align*}
&\sigma^2(x) - \sigma^2(c) 
\\
&\le 
\Bigl(x(1-x) + \frac{2}{3\cdot 4} (3x-1)(2-4x)
+ \frac{2}{3^2\cdot 4^2} (4^2x-6)(4-3^2x)
\Bigr)\Bigr|_{x={27}/{64}}
 + \frac1{22\cdot 12^2}- \sigma^2(c) 
\\&
=
-\frac{1477989115}{5393688588288}
<0
,
\end{align*}
since the quadratic function has the axis  
at $x={31}/{72} > {27}/{64}$.

On $[\,{16}/{37}, {7}/{16})$,  by  noting $c < 16/37$ 
and by applying (\ref{Eq:est}) for $N=2$ we have
\begin{align*}
&\sigma^2(x) - \sigma^2(c) 
\\
&\le 
\Bigl(x(1-x) + \frac{2}{3\cdot 4} (3x-1)(2-4x)
+ \frac{2}{3^2\cdot 4^2} (3^2x-3)(7-4^2x)
\Bigr)\Bigr|_{x={16}/{37}}
 + \frac1{22\cdot 12^2}- \sigma^2(c) 
\\&
=
-\frac{16743731111}{57687184979424}
<0
,
\end{align*}
since the quadratic function has the axis  
at $x={101}/{240} < 16/37$.

\subsection{ $[\,{27}/{64}, c)$ and $[\,c, {16}/{37})$  parts.}
On $[\,{27}/{64}, {16}/{37})$, 
we have
$\langle 4^3 x\rangle = 4^3 x -27$, 
$\langle 3^2 x\rangle = 3^3 x -11$, and
$\langle 4^3 x\rangle - \langle 3^3 x\rangle = 37x -16 < 0$. 

On $[\,{27}/{64},c)$, 
by applying (\ref{Eq:dest}) for $N=3$ 
we have 
\begin{align*}
\frac{d}{dx} \sigma^2(x) 
&\ge 
\frac{d}{dx} 
\Bigl(
x(1-x) + \frac{2}{3\cdot 4} (3x-1)(2-4x)
+ \frac{2}{3^2\cdot 4^2} (4^2x-6)(4-3^2x)
\\&\qquad\qquad
+ \frac{2}{3^3\cdot 4^3} (4^3x-27)(12-3^3x)
\Bigr)-\frac{1}{3^3}
\\&=
\Bigl(\frac{1739}{288} - 14 x\Bigr) -\frac{1}{3^3}
> \Bigl(\frac{1739}{288} - 14 c\Bigr) -\frac{1}{3^3}
=\frac{1}{864}
>0.
\end{align*}

On $[\,c, {16}/{37})$,
by applying (\ref{Eq:dest}) for $N=3$ 
we have 
\begin{align*}
\frac{d}{dx} \sigma^2(x) 
&\le 
\frac{d}{dx} 
\Bigl(
x(1-x) + \frac{2}{3\cdot 4} (3x-1)(2-4x)
+ \frac{2}{3^2\cdot 4^2} (3^2x-3)(7-4^2x)
\\&\qquad\qquad
+ \frac{2}{3^3\cdot 4^3} (4^3x-27)(12-3^3x)
\Bigr)+\frac{1}{3^3}
\\&=
\Bigl(\frac{1711}{288} - 14 x\Bigr) +\frac{1}{3^3}
< \Bigl(\frac{1711}{288} - 14 c\Bigr) +\frac{1}{3^3}
=-\frac{19}{864}
<0
.
\end{align*}

\section{Type II case, $\Sigma_{8/3}$}

Put
$
c = \dfrac{24}{55}
$. 
By $3^{20}=8^{20}=1$ mod 55, 
we see
$V\bigl(\langle 3^{k+20j}c\rangle, \langle 8^{k+20j}c\rangle\bigr) 
=V\bigl(\langle 3^{k}c\rangle, \langle 8^{k}c\rangle\bigr) 
$
and
\begin{align*}
\sigma^2(c)&=
V(c, c)
+
\frac{2(3\cdot 8)^{20}}{(3\cdot 8)^{20}-1}\sum_{k=1}^{20}\frac1{(3\cdot8)^{k}} 
V\bigl(\langle 3^{k}c\rangle, \langle 10^{k}c\rangle\bigr) 
\\&=
\frac{630671157052050732205279776888
}{2432093174293565222595046931875
}
.
\end{align*}
We divide $[\,0,1/2)$ into 
$[\,0,\afrac38)$, 
$[\,\afrac38,\afrac25)$, 
$[\,\afrac25, \afrac{27}{8^2})$, 
$[\,\afrac{27}{8^2},\afrac{223}{8^3})$, 
$[\,\afrac{223}{8^3}, c)$, 
$[\,c, \afrac{28}{8^2})$, 
$[\,\afrac{28}{8^2}, \afrac{29}{8^2})$, 
$[\,\afrac{29}{8^2}, \afrac{30}{8^2})$, 
$[\,\afrac{30}{8^2}, \afrac12)$, 
and prove $ \sigma^2(x) <  \sigma^2(c) $ ($x\not=c$) on each. 

\subsection{ $[\,0,\afrac38)$ part.}
By applying (\ref{Eq:est}) for $N=0$ we have
\begin{align*}
 \sigma^2(x) -  \sigma^2(c) 
&\le 
x(1-x) + \frac{1}{46} -  \sigma^2(c) 
\le
\Bigl(x(1-x) + \frac{1}{46} -  \sigma^2(c) \Bigr)\Bigm|_{x=3/8}
\\
&=
-\frac{497774629214982125558223402707}{155653963154788174246083003640000}
<0.
\end{align*}

\subsection{ $[\,\afrac38,\afrac25)$, $[\,\afrac25, \afrac{27}{8^2})$ and 
 $[\,\afrac{30}{8^2}, \afrac12)$ parts.}

On $[\,\afrac38, \afrac12)$, 
we have
$\langle 8 x\rangle = 8x-3$, 
$\langle 3 x\rangle = 3x-1$, 
$\langle 8 x\rangle -\langle 3 x\rangle = 5x-2<0$ if and only if $x< \afrac25$, and
\begin{equation}\label{Eq:3/8second}
V\bigl(\langle 3 x\rangle, \langle 8 x\rangle\bigr)
=
\begin{cases}
(8x-3)(2-3x) & x\in [\,3/8,2/5),
\\
(3x-1)(4-8x) & x\in [\,2/5,1/2).
\end{cases}
\end{equation}

On $[\,\afrac38, \afrac25)$, 
by applying (\ref{Eq:est}) for $N=1$ we have
\begin{align*}
 \sigma^2(x) -  \sigma^2(c) 
&\le 
x(1-x) +  \frac{2}{3\cdot 8}(8x-3)(2-3x) + \frac{1}{46\cdot 24}-  \sigma^2(c) 
\\&
\le
\Bigl(x(1-x) +  \frac{2}{3\cdot 8}(8x-3)(2-3x) + \frac{1}{46\cdot 24}-  \sigma^2(c) \Bigr)
\Bigm|_{x=\afrac25}
\\
&=
-\frac{197406452532552488067407258633}
{38913490788697043561520750910000}
<0, 
\end{align*}
since the quadratic function has the axis at $x=\afrac{37}{72}> \afrac25$.

On $[\,\afrac25, \afrac{27}{8^2})$, 
by applying (\ref{Eq:est}) for $N=1$  
we have
\begin{align*}
 \sigma^2(x) -  \sigma^2(c) 
&\le 
\Bigl(x(1-x) +  \frac{2}{3\cdot 8}(3x-1)(4-8x) + \frac{1}{46\cdot 24}-  \sigma^2(c)\Bigr)
\Bigm|_{x=\afrac{27}{8^2}}
\\
&=
-\frac{20177840406176046206732438431619}{29885560925719329455247936698880000}
<0, 
\end{align*}
since the  quadratic function has the axis at $x=\afrac49\in (\afrac{27}{8^2},\afrac{30}{8^2})$. 
On $[\,\afrac{30}{8^2},\afrac12)$, it also implies
\begin{align*}
 \sigma^2(x) -  \sigma^2(c) 
&\le 
\Bigl(x(1-x) +  \frac{2}{3\cdot 8}(3x-1)(4-8x) + \frac{1}{46\cdot 24}-  \sigma^2(c)\Bigr)
\Bigm|_{x=\afrac{30}{8^2}}
\\
&=
-\frac{6868529982264185468629394806811}{7471390231429832363811984174720000}
<0.
\end{align*}

\subsection{ $[\,\afrac{27}{8^2},\afrac{223}{8^3})\subset[\,\afrac{27}{8^2},c)$ part.}
We have
$\langle 8^2 x\rangle = 8^2x -27$, 
$\langle 3^2 x\rangle = 3^2x -3$, and
$\langle 8^2 x\rangle -\langle 3^2 x\rangle =55x-24< 0$.
By applying (\ref{Eq:est}) for $N=2$ we  have
\begin{align*}
& \sigma^2(x) -  \sigma^2(c) 
\\
&\le 
\Bigl(x(1-x) +  \frac{2}{3\cdot 8}(3x-1)(4-8x) +
\frac{2}{3^2\cdot 8^2}(8^2x-27)(4-3^2x)
+ \frac{1}{46\cdot 24^2}-  \sigma^2(c)\Bigr)
\Bigm|_{x=\afrac{223}{8^3}}
\\
&=
-\frac{53326981170609906995809448895223
}{5738027697738111255407603846184960000}
<0,  
\end{align*}
since the quadratic function has the axis at $x=\afrac{1267}{2880}> \afrac{223}{8^3}$. 

\subsection{ $[\,\afrac{29}{8^2}, \afrac{30}{8^2})$ part.}
We have
$\langle 8^2 x\rangle = 8^2x -29$, 
$\langle 3^2 x\rangle = 3^2x -4$, and 
$\langle 8^2 x\rangle -\langle 3^2 x\rangle =55x-25<  0$ if and only if  
$x< 5/11\in (\afrac{29}{8^2}, \afrac{30}{8^2})$. 

On $[\,\afrac{29}{8^2}, \afrac{5}{11})$,
by applying (\ref{Eq:est}) for $N=2$  we have
\begin{align*}
&\sigma^2(x) -  \sigma^2(c) 
\\
&\le 
\Bigl(x(1-x) +  \frac{2}{3\cdot 8}(3x-1)(4-8x) +
\frac{2}{3^2\cdot 8^2}(8^2x-29)(4-3^2x)
+ \frac{1}{46\cdot 24^2}-  \sigma^2(c)\Bigr)
\Bigm|_{x=\afrac{5}{11}}
\\
&=
-\frac{10652022559016351057295568539
}{311307926309576348492166007280000}
<0,  
\end{align*}
since the quadratic function has the axis at $x=\afrac{1349}{2880}$.

On $[\, \afrac{5}{11},\afrac{30}{8^2})$,
by applying (\ref{Eq:est}) for $N=2$  we have
\begin{align*}
& \sigma^2(x) -  \sigma^2(c) 
\\
&\le 
\Bigl(x(1-x) +  \frac{2}{3\cdot 8}(3x-1)(4-8x) +
\frac{2}{3^2\cdot 8^2}(3^2x-4)(30-8^2x)
+ \frac{1}{46\cdot 24^2}-  \sigma^2(c)\Bigr)
\Bigm|_{x=\afrac{5}{11}}
\\
&=
\frac{10652022559016351057295568539
}{311307926309576348492166007280000}
<0,  
\end{align*}
since the quadratic function has the axis at $x=\afrac{647}{1440}<\afrac{5}{11}$.

\subsection{ $[\,\afrac{28}{8^2}, \afrac{29}{8^2})$ part.}
On $[\,\afrac{28}{8^2}, \afrac{4}{3^2})$, 
we have $\langle 8^2 x\rangle = 8^2x -28$, 
$\langle 3^2 x\rangle = 3^2x -3$, and
$\langle 8^2 x\rangle -\langle 3^2 x\rangle = 55x-25<0$.
By applying (\ref{Eq:est}) for $N=2$ we have
\begin{align*}
& \sigma^2(x) -  \sigma^2(c) 
\\
&\le 
\Bigl(x(1-x) +  \frac{2}{3\cdot 8}(3x-1)(4-8x) +
\frac{2}{3^2\cdot 8^2}(8^2x-28)(4-3^2x)
+ \frac{1}{46\cdot 24^2}-  \sigma^2(c)\Bigr)
\Bigm|_{x=\afrac{319}{720}}
\\
&=
-\frac{274411727420190945072275586943
}{50431884062151368455730893179360000
}
<0,  
\end{align*}
since the quadratic function has the axis at $x=\afrac{319}{720}$.

On $[\,\afrac{4}{3^2}, \afrac{29}{8^2} )$, 
we have $\langle 8^2 x\rangle = 8^2x -28$, 
$\langle 3^2 x\rangle = 3^2x -4$, and 
$\langle 8^2 x\rangle -\langle 3^2 x\rangle = 55x-24>0$ by $\afrac{24}{55}<  \afrac{4}{3^2}$. 
By applying (\ref{Eq:est}) for $N=2$ we have
\begin{align*}
& \sigma^2(x) -  \sigma^2(c) 
\\
&\le 
\Bigl(x(1-x) +  \frac{2}{3\cdot 8}(3x-1)(4-8x) +
\frac{2}{3^2\cdot 8^2}(3^2x-4)(29-8^2x)
+ \frac{1}{46\cdot 24^2}-  \sigma^2(c)\Bigr)
\Bigm|_{x=\afrac{257}{576}}
\\
&=
-\frac{12819924994637720485324913713
}{806910144994421895291694290869760000}
<0,  
\end{align*}
since the quadratic function has the axis at $x=\afrac{257}{576}$.

\subsection{ $[\,\afrac{223}{8^3}, c)$ and $[\,c, \afrac{28}{8^2})$ parts.}

We have
$\langle 8^2 x\rangle = 8^2x -27$, 
$\langle 3^2 x\rangle = 3^2x -3$, and 
$\langle 8^2 x\rangle -\langle 3^2 x\rangle = 55x -24 <0$ if and only if $x< c$. 

On $[\,c, \afrac{28}{8^2})$, 
by applying (\ref{Eq:dest}) for $N=2$  we have
\begin{align*}
\frac{d}{dx} \sigma^2(x) 
&\le 
\frac{d}{dx} 
\Bigl(
x(1-x) + \frac{2}{3\cdot 8} (3x-1)(4-8x) 
+ \frac{2}{3^2\cdot 8^2} (3^2x-3)(28-8^2x)
\Bigr)+\frac{1}{3^2}
\\&=
\Bigl(\frac{311}{72} - 10 x\Bigr) 
\le 
\Bigl(\frac{311}{72} - 10 c\Bigr) 
=-\frac{35}{792}
<0. 
\end{align*}

On $[\,\afrac{223}{8^3},c)$, 
we have
$\langle 8^3 x\rangle = 8^3x -223$, 
$\langle 3^3 x\rangle = 3^3x -11$, and 
$\langle 8^3 x\rangle -\langle 3^3 x\rangle = 485x -212 <0$.  
By applying (\ref{Eq:dest}) for $N=3$ we have
\begin{align*}
\frac{d}{dx} \sigma^2(x) 
&\ge 
\frac{d}{dx} 
\Bigl(
x(1-x) + \frac{2}{3\cdot 8} (3x-1)(4-8x) 
+ \frac{2}{3^2\cdot 8^2} (8^2x-27)(4-3^2x)
\\&\qquad\qquad\qquad
+ \frac{2}{3^3\cdot 8^3} (8^3x-223)(12-3^2x)
\Bigr)-\frac{1}{3^3}
\\&=
\Bigl(\frac{42317}{6912} - 14 x\Bigr) 
\ge 
\Bigl(\frac{42317}{6912} - 14 c\Bigr) 
=\frac{5003}{380160}
>0. 
\end{align*}

\section{Type II case, $\Sigma_{10/3}$}

Put
$
c = \dfrac{40}{91}
$. 
By $3^{6}=10^{6}=1$ mod 91, 
we see
$V\bigl(\langle 3^{k+6j}c\rangle, \langle 10^{k+6j}c\rangle\bigr) 
=V\bigl(\langle 3^{k}c\rangle, \langle 10^{k}c\rangle\bigr) 
$
and
\begin{align*}
\sigma^2(c)&=
V(c, c)
+
\frac{2(3\cdot 10)^{6}}{(3\cdot 10)^{6}-1}\sum_{k=1}^{6}\frac1{(3\cdot10)^{k}} 
V\bigl(\langle 3^{k}c\rangle, \langle 10^{k}c\rangle\bigr) 
%\\&
=
\frac{1565277378120}{6036848991719}
.
\end{align*}

We divide $[\,0,1/2)$ into 
$[\,0, {4}/{10})$, 
$[\,4/{10}, {43}/{10^2}\,)$, 
$[\,{43}/{10^2}, {4389}/{10^4})$,
$[\,{4389}/{10^4}, {439}/{10^3})$, \break
$[\,{439}/{10^3}, {4395}/{10^4})$,
$[\,{4395}/{10^4}, c)$,
$[\,c, {44}/{10^2})$,
$[\, {44}/{10^2}, {4}/{3^2})$,
$[\,{4}/{3^2}, {45}/{10^2})$,
$[\,{45}/{10^2}, {46}/{10^2})$, 
$[\,{46}/{10^2}, 1/2)$,
and prove $\sigma^2(x) < \sigma^2(c) $ ($x\not=c$) on each. 

\subsection{ $[\,0, {4}/{10})$ part.}
By applying (\ref{Eq:est}) for $N=0$ we have
\begin{align*}
\sigma^2(x) - \sigma^2(c) 
\le  x(1-x) +\frac1{58} - \sigma^2(c) \biggl|_{x= {4}/{10}}
=
-\frac{617500840097}{301842449585950}                                  
< 0.
\end{align*}

\subsection{ $[\,4/{10}, {43}/{10^2}\,)$ and  $[\,{46}/{10^2}, 1/2)$ parts.}

On $[\,4/{10}, 1/{2}\,)$, 
we have
$\langle 10x\rangle= 10x -4$, $\langle 3x \rangle = 3x-1$, 
$\langle 10x\rangle- \langle 3x \rangle = 7x -3$, and hence
\begin{equation}\label{Eq:10/3second}
V(\langle 3x \rangle, \langle 10x\rangle)
=
\begin{cases}
(10x-4)(2-3x) & x\in  [\, 4/{10}, 3/7),
\\
(3x-1)(5-10x) & x\in  [\, 3/7, 1/2).
\end{cases}
\end{equation}

On $[\, 4/{10}, 3/7)$, 
by applying (\ref{Eq:est}) for $N=1$ 
we have
\begin{align*}
\sigma^2(x) - \sigma^2(c) 
&\le x(1-x)  + \frac{2}{3\cdot 10}(10x-4)(2-3x) 
+\frac1{58\cdot 30} - \sigma^2(c) 
\biggl|_{x=3/7}
\\
&=
-\frac{25238536343}{120736979834380}
< 0,
\end{align*}
since
the quadratic function has the axis at 
$x={47}/{90}> 3/7$.

On $[\,3/7, {43}/{10^2})$, by applying (\ref{Eq:est}) for $N=1$ we have
\begin{align*}
\sigma^2(x) - \sigma^2(c) 
&\le x(1-x)  + \frac{2}{3\cdot 10}(3x-1)(5-10x)
+\frac1{58\cdot 30} - \sigma^2(c) 
\biggl|_{x={43}/{10^2}}
\\&
=
-\frac{14326413546779}{181105469751570000}< 0
,
\end{align*}
since the quadratic function has the axis at
$x=4/9\in ({43}/{10^2}, {46}/{10^2})$. 
On $[\,{46}/{10^2}, 1/2)$, it also implies 
\begin{align*}
\sigma^2(x) - \sigma^2(c) 
&\le x(1-x)  + \frac{2}{3\cdot 10}(3x-1)(5-10x)
+\frac1{58\cdot 30} - \sigma^2(c) 
\biggl|_{x={46}/{10^2}}
\\&
=
-\frac{2027310032621}{11319091859473125}
< 0.
\end{align*}

\subsection{ $[\,{43}/{10^2}, {4389}/{10^4})$ part.}

On $[\,{43}/{10^2}, c)$, we have
$\langle 3^2 x\rangle = 3^2x-3$, 
$\langle 10^2 x\rangle = 10^2x-43$, 
$\langle 10^2 x\rangle - \langle 3^2 x\rangle = 91 x -40< 0$, and
\begin{equation}\label{Eq:10/3third}
V(\langle 3^2 x\rangle, \langle 10^2 x\rangle)
=(10^2x -43)(4-3^2x) 
\quad
x\in [\,{43}/{10^2}, c).
\end{equation}
On $[\,{43}/{10^2}, {4389}/{10^4})\subset[\,{43}/{10^2}, c)$, 
by applying (\ref{Eq:est}) for $N=2$  we have
\begin{align*}
\sigma^2(x) - \sigma^2(c) 
&\le x(1-x)  + \frac{2}{3\cdot 10}(3x-1)(5-10x)
+\frac{2}{3^2\cdot 10^2}(10^2x-43)(4-3^2x)
\\
&\qquad
+\frac1{58\cdot 30^2} - \sigma^2(c) 
\biggr|_{x={4389}/{10^4}}
\\
&=
-\frac{2462822876208511}{1086632818509420000000}< 0
,
\end{align*}
since the quadratic function has the axis  at 
$x={1987}/{4500}> {4389}/{10^4}$.

\subsection{ $[\,{4389}/{10^4}, {439}/{10^3})$ part.}
On $[\,{4389}/{10^4}, {439}/{10^3})$, we have
$\langle 3^3 x\rangle = 3^3x-11$, 
$\langle 10^3 x\rangle = 10^3x-438$, and
$\langle 10^3 x\rangle - \langle 3^3 x\rangle = 973 x -427>  0$. 
By applying (\ref{Eq:est}) for $N=3$ we have
\begin{align*}
\sigma^2(x) - \sigma^2(c) 
&\le x(1-x)  + \frac{2}{3\cdot 10}(3x-1)(5-10x)
+\frac{2}{3^2\cdot 10^2}(10^2x-43)(4-3^2x)
\\
&\qquad
+\frac{2}{3^3\cdot 10^3}(3^3x-11)(439-10^3x)
+\frac1{58\cdot 30^3} - \sigma^2(c) 
\biggr|_{x={4389}/{10^4}}
\\&=
-\frac{78707379591968117}{5433164092547100000000}< 0
,
\end{align*}
since
the quadratic function has the axis at 
$x={82463}/{189000} < {4389}/{10^4}$.

\subsection{ $[\,{439}/{10^3}, {4395}/{10^4})$ part.}
We have
$\langle 3^3 x\rangle = 3^3x-11$, 
$\langle 10^3 x\rangle = 10^3x-439$, and
$\langle 10^3 x\rangle - \langle 3^3 x\rangle = 973 x -428<   0$. 
By applying (\ref{Eq:est}) for $N=3$ we have
\begin{align*}
\sigma^2(x) - \sigma^2(c) 
&\le x(1-x)  + \frac{2}{3\cdot 10}(3x-1)(5-10x)
+\frac{2}{3^2\cdot 10^2}(10^2x-43)(4-3^2x)
\\
&\qquad
+\frac{2}{3^3\cdot 10^3} (10^3x-439)(12-3^3x)
+\frac1{58\cdot 30^3} - \sigma^2(c) 
\biggr|_{x={4395}/{10^4}}
\\&=
-\frac{1113068046625279}{651979691105652000000}< 0
,
\end{align*}
since
the quadratic function has the axis at
    $x={27821}/{63000} > {4395}/{10^4}$.

\subsection{ $[\,{45}/{10^2}, {46}/{10^2})$ part.}
We have
$\langle 3^2 x\rangle = 3^2x-4$, 
$\langle 10^2 x\rangle = 10^2x-45$, and
$\langle 10^2 x\rangle - \langle 3^2 x\rangle = 91 x -41$. 
Note that ${41}/{91}$ is in this interval. 
Therefore 
$$
V(\langle 3^2 x\rangle , \langle 10^2 x\rangle )
=
\begin{cases}
(10^2 x-45)(5-3^2x) & x\in [\,{45}/{10^2}, {41}/{91}), 
\\
(3^2x-4)(46-10^2x) & x\in [\,{41}/{91}, {46}/{10^2}). 
\end{cases}
$$
On $[\,{41}/{91}, {46}/{10^2})$, 
by applying (\ref{Eq:est}) for $N=2$ 
we have
\begin{align*}
\sigma^2(x) - \sigma^2(c) 
&\le x(1-x)  + \frac{2}{3\cdot 10}(3x-1)(5-10x)
+\frac{2}{3^2\cdot 10^2}(3^2x-4)(46-10^2x)
\\
&\qquad
+\frac1{58\cdot 30^2} - \sigma^2(c) 
\biggr|_{x={41}/{91}}
\\&=
-\frac{18679091303
}{3622109395031400
}< 0
,\end{align*}
since
the quadratic function has the axis at
    $x={1007}/{2250}< {41}/{91}$.

On $[\, {45}/{10^2}, {41}/{91})$, 
by applying (\ref{Eq:est}) for $N=2$ 
we have
\begin{align*}
\sigma^2(x) - \sigma^2(c) 
&\le x(1-x)  + \frac{2}{3\cdot 10}(3x-1)(5-10x)
+\frac{2}{3^2\cdot 10^2}(10^2x-45)(5-3^2x)
\\
&\qquad
+\frac1{58\cdot 30^2} - \sigma^2(c) 
\biggr|_{x={41}/{91}}
\\&=
-\frac{18679091303
}{3622109395031400
}< 0
,\end{align*}
since the quadratic function has the axis at
    $x={421}/{900}
> {41}/{91}$.

\subsection{ $[\,{4}/{3^2}, {45}/{10^2})$ part.}
We see
$\langle 3^2 x\rangle = 3^2x-4$, 
$\langle 10^2 x\rangle = 10^2x-44$, and
$\langle 10^2 x\rangle - \langle 3^2 x\rangle = 91 x -40> 0$. 
By applying (\ref{Eq:est}) for $N=2$ we have
\begin{align*}
\sigma^2(x) - \sigma^2(c) 
&\le x(1-x)  + \frac{2}{3\cdot 10}(3x-1)(5-10x)
+\frac{2}{3^2\cdot 10^2}(3^2x-4)(45-10^2x)
\\
&\qquad
+\frac1{58\cdot 30^2} - \sigma^2(c) 
\biggr|_{x={401}/{900}}
\\&=
-\frac{2505705655291}{977969536658478000}< 0
,\end{align*}
since
the quadratic function has the axis at
    $x={401}/{900}$. 

\subsection{ $[\, {44}/{10^2}, {4}/{3^2})$ part.}

We see
$\langle 3^2 x\rangle = 3^2x-3$, 
$\langle 10^2 x\rangle = 10^2x-44$, and
$\langle 10^2 x\rangle - \langle 3^2 x\rangle = 91 x -41< 0$.
By applying (\ref{Eq:est}) for $N=2$ we have
\begin{align*}
\sigma^2(x) - \sigma^2(c) 
&\le x(1-x)  + \frac{2}{3\cdot 10}(3x-1)(5-10x)
+\frac{2}{3^2\cdot 10^2}(10^2x-44)(4-3^2x)
\\
&\qquad
+\frac1{58\cdot 30^2} - \sigma^2(c) 
\biggr|_{x={499}/{1125}}
\\&=
-\frac{58487141153873}{12224619208230975000}< 0
,\end{align*}
since
the quadratic function has the axis at
    $x={499}/{1125}$. 

\subsection{ $[\,{4395}/{10^4}, c)$ and $[\,c, {44}/{10^2})$  parts.}
We have
$\langle 3^2 x\rangle = 3^2x -3$, 
$\langle 10^2 x\rangle = 10^2x -43$, hence  
$\langle 10^2 x\rangle - \langle 3^2 x\rangle = 91x -40 > 0 $
if and only if $x> c$. 
Therefore 
$$
V\bigl(\langle 3^2 x\rangle, \langle 10^2 x\rangle\bigr)
=
\begin{cases}
(10^2x-43)(4-3^2x) & x\in [\,{4395}/{10^4}, c),
\\
(3^2x-3)(44-10^2x) & x\in [\,c, {44}/{10^2}).
\end{cases}
$$
On $[\,c,{44}/{10^2})$, 
by applying (\ref{Eq:dest}) for $N=2$ we have
\begin{align*}
\frac{d}{dx} \sigma^2(x) 
&\le 
\frac{d}{dx} 
\Bigl(
x(1-x) +\frac{2}{3\cdot 10} (3x-1)(5-10x) 
+ \frac{2}{3^2\cdot 10^2} (3^2x-3)(44-10^2x)
\Bigr)+\frac{1}{3^2}
\\&=
\Bigl(\frac{973}{225} - 10 x\Bigr) 
\ge 
\Bigl(\frac{973}{225} - 10 c\Bigr) 
=-\frac{1457}{20475}
<0. 
\end{align*}

On $[\,{4395}/{10^4},c)$
we have
$\langle 3^3 x\rangle = 3^3x -11$, 
$\langle 10^3 x\rangle = 10^3x -439$, 
$\langle 10^3 x\rangle - \langle 3^3 x\rangle = 973x -428 < 0 $, 
$\langle 3^4 x\rangle = 3^4x -35$, 
$\langle 10^4 x\rangle = 10^4x -4395$, and
$\langle 10^4 x\rangle - \langle 3^4 x\rangle = 9919x -4360 < 0 $
since $\afrac{4360}{9919}= c$.
By applying (\ref{Eq:dest}) for $N=4$, we have
\begin{align*}
\frac{d}{dx} \sigma^2(x) 
&\ge 
\frac{d}{dx} 
\Bigl(
x(1-x) +\frac{2}{3\cdot 10} (3x-1)(5-10x) 
+ \frac{2}{3^2\cdot 10^2} (10^2x-43)(4-3^2x)
\\&\qquad
+\frac{2}{3^3\cdot 10^3} (10^3x -439)(12-3^3x)
+\frac{2}{3^4\cdot 10^4} (10^4x -4395)(36-3^4x)
\Bigr)-\frac{1}{3^4}
\\&=
\frac{642977}{81000}-18x
\ge 
\frac{642977}{81000}-18c
=
\frac{190907}{7371000}
>0. 
\end{align*}

\section{Type II case, $\Sigma_{12/5}$}

Put $
c = \dfrac{55}{119}$.
Since we have $12^{48}=1$, $5^{48}=1$ mod $119$, 
we have
$\langle 12^{48k+n}c\rangle = \langle 12^{n}c\rangle$, 
$\langle 5^{48k+n}c\rangle = \langle 5^{n}c\rangle$, 
and
\begin{align*}
&\sigma^2(c)
=
V(c,c) 
+
2\frac{12^{48}5^{48}}{12^{48}5^{48}-1}
\sum_{n=1}^{48}
\frac1{12^n5^n}
V\bigl(\langle 12^{n}c\rangle, \langle 5^{n}c\rangle\bigr)
\\
&=
\tfrac{80571343410682351163664747839964625093429559742738931998565277817228239636988823752192520}
{317946421393847885076874996724028068069375999999999999999999999999999999999999999999985839}
.
\end{align*}

We divide $[\,0,1/2)$ into 
$[0, 3/7)$, 
$[\,3/7,66/12^2)$, 
$[\,66/12^2, c)$, 
$[\,c, 67/12^2)$, 
$[\,67/12^2, 68/12^2)$, 
$[\,68/12^2,1/2)$, 
and prove $\sigma^2(x) < \sigma^2(c) $ ($x\not=c$) on each.

\subsection{$[0, 3/7)$ part.}

By applying (\ref{Eq:est}) for $N=0$  we have
\begin{align*}
&\sigma^2(x) -\sigma^2(c) 
< 
x(1-x) +\frac{1}{2\cdot 59}-\sigma^2(c) \biggm|_{x=3/7}
\\&
=
 - 
\tfrac{24904966647596952849429889302781933093695756664304675096657329371733511265783097612555}
{635892842787695770153749993448056136138751999999999999999999999999999999999999999999971678}
\\&<0.
\end{align*}

\subsection{ $[\,3/7,66/12^2)$ and  $[\,68/12^2,1/2)$ parts.} 

On $[\,3/7, 1/2)$, 
we have $\langle 5 x\rangle= 5x-2 < 12x-5= \langle 12 x\rangle$ and 
\begin{equation}\label{Eq:12/5second}
V\bigl(\langle 5 x\rangle, \langle 12 x\rangle\bigr)
=
(5x-2)(6-12x), 
\quad
x\in [\,3/7, 1/2).
\end{equation}

On $[\,3/7,66/12^2)$,  by applying (\ref{Eq:est}) for $N=1$  
we have 
\begin{align*}
&\sigma^2(x) -\sigma^2(c) \le
x(1-x) +\frac{2}{5\cdot 12}(5x-2)(6-12x)
+\frac{1}{2\cdot 59\cdot 5\cdot 12}\biggm|_{x=66/12^2}
\\&
=
-\tfrac{8879579275594774515930473387843148598818792639773248809279103619686417081515177370804991}
{61045712907618793934759999371013389069320191999999999999999999999999999999999999999997281088}
< 0, 
\end{align*}
since the quadratic function has the axis at $x=7/15 \in [\,66/12^2, 68/12^2)$. 
On $[\,68/12^2,1/2)$, it also implies 
\begin{align*}
&\sigma^2(x) -\sigma^2(c) \le
x(1-x) +\frac{2}{5\cdot 12}(5x-2)(6-12x)
+\frac{1}{2\cdot 59\cdot 5\cdot 12}\biggm|_{x=68/12^2}
\\&=
-\tfrac{20408661501979242034999076432228404719367417197449049104389915721472192167045745421559689}
{686764270210711431766049992923900627029852159999999999999999999999999999999999999999969412240}
< 0.
\end{align*}

\subsection{ $[\,67/12^2, 68/12^2)$ part.}
On $[\,67/12^2, 68/12^2)$ we have
$\langle 12^2 x\rangle=12^2x-67$, 
$\langle 5^2 x\rangle=5^2 x-11$, 
$\langle 12^2 x\rangle-\langle 5^2 x\rangle=119x-56$, 
$56/119\in [\,67/12^2, 68/12^2)$, and
$$
V\bigl(\langle 5^2 x\rangle, \langle 12^2 x\rangle\bigr)
=
V\bigl(
=
\begin{cases}
(12^2x-67)(12-5^2x) & x \in [\,67/12^2, 56/119),
\\
(5^2 x-11)(68-12^2x)& x\in [\,56/119,68/12^2).
\end{cases}
$$

On $ [\,67/12^2, 56/119)$,
 by applying (\ref{Eq:est}) for $N=2$   we have
\begin{align*}
&\sigma^2(x) -\sigma^2(c) \le 
x(1-x) +\frac{2}{5\cdot 12}(5x-2)(6-12x)
+\frac{2}{5^2 12^2}(12^2x-67)(12-5^2x)
\\
&\qquad\qquad\qquad\qquad
+\frac{1}{2\cdot 59\cdot 5^2 12^2}
\biggm|_{x=8443/18000}-\sigma^2(c) 
\\
&
=
-
\tfrac{215144764738461252588955422058541084562371891923471473131697471644165765011372362646808357089}
{20602928106321342952981499787717018810895564799999999999999999999999999999999999999999082367200000}
\\&
<0,
\end{align*}
since the quadratic function has the axis at 
$x=8443/18000$.

On $ [\,56/119,68/12^2)$, 
 by applying (\ref{Eq:est}) for $N=2$   we have
\begin{align*}
&\sigma^2(x) -\sigma^2(c) \le 
x(1-x) +\frac{2}{5\cdot 12}(5x-2)(6-12x)
+\frac{2}{5^2 12^2}(5^2 x-11)(68-12^2x)
\\
&\qquad\qquad\qquad\qquad
+\frac{1}{2\cdot 59\cdot 5^2 12^2}-\sigma^2(c) 
\biggm|_{x=56/119}
\\
&
=
-
\tfrac{16931007279972763891785546344909946343602907997165610115988795246080213518939717135066577}
{763071411345234924184499992137667363366502399999999999999999999999999999999999999999966013600}
\\&
<0,
\end{align*}
since the  quadratic function has the axis at 
$x=2081/4500<56/119$.

\subsection{ $[\,66/12^2, c)$ and $[\,c, 67/12^2)$ parts.}
We have
$\langle 12^2 x\rangle=12^2x-66$, 
$\langle 5^2 x\rangle=5^2 x-11$, 
$\langle 12^2 x\rangle-\langle 5^2 x\rangle=119x-55$, 
$c\in [\,66/12^2, 67/12^2)$, and 
\begin{equation}\label{Eq:12/5third}
V\bigl(\langle 5^2 x\rangle, \langle 12^2 x\rangle\bigr)
=
\begin{cases}
(12^2x-66)(12-5^2x) & x \in [\,66/12^2, c),
\\
(5^2 x-11)(67-12^2x)& x\in [\,c,67/12^2),
\end{cases}
\end{equation}
We use (\ref{Eq:dest}) for $N=2$. 

On $[\,66/12^2,c)$, 
we have 
\begin{align*}
\frac{d}{dx}
\sigma^2(x)
&\ge
\frac{d}{dx}
\Bigl(
x(1-x)
+
\frac{2}{5\cdot 12}(5x-2)(6-12x)
+
\frac{2}{5^2 12^2}(12^2x-66)(12-5^2x)
\Bigr)
-\frac{1}{2\cdot 5^2}
\\
&
=
-\frac{3000 x - 1397}{300}
\ge 
-\frac{3000c - 1397}{300}
=
\frac{1243}{35700}
> 0
.
\end{align*}

Note that $801/12^3\in [\,c, 67/12^2)$. 
On $[\,801/12^3, 67/12^2)$,  
 by applying (\ref{Eq:dest}) for $N=2$  we have 
\begin{align*}
\frac{d}{dx}
\sigma^2(x)
&\le
\frac{d}{dx}
\Bigl(
x(1-x)
+
\frac{2}{5\cdot 12}(5x-2)(6-12x)
+
\frac{2}{5^2 12^2}(5^2 x-11)(67-12^2x)
\Bigr)
+\frac{1}{2\cdot 5^2}
\\
&
=
-\frac{3600 x - 1667}{360}
\ge 
-\frac{3600 (801/12^3) - 1667}{360}
=
-\frac{7}{1440}
<0
. 
\end{align*}

Now and later on we apply  (\ref{Eq:dest}) for $N=3$. 

On $[\, 798/12^3, 799/12^3)$, we have
$\langle 12^3 x\rangle=12^3x-798$, 
$\langle 5^3 x\rangle=5^3x-57$, 
$\langle 12^3 x\rangle-\langle 5^3 x\rangle=1603x-741$, 
$741/1603 \in [\,c, 799/12^3)$, and
$$
V\bigl(\langle 5^3 x\rangle, \langle 12^3 x\rangle\bigr)
=
\begin{cases}
(12^3x-798)(58-5^3x) & x\in [\, 798/12^3, 741/1603),
\\
(5^3x-57)(799-12^3x) & x\in [\,  741/1603, 799/12^3).
\end{cases}
$$
On $[\, c, 741/1603)$ 
we have 
\begin{align*}
\frac{d}{dx}
\sigma^2(x)
&\le
\frac{d}{dx}
\Bigl(
x(1-x)
+
\frac{2}{5\cdot 12}(5x-2)(6-12x)
+
\frac{2}{5^2 12^2}(5^2 x-11)(67-12^2x)
\\&\qquad
+
\frac{2}{5^3 12^3}(12^3x-798)(58-5^3x)
\Bigr)
+\frac{1}{2\cdot 5^3}
\\
&
=
-\frac{84000 x - 38797}{6000}
\le 
-\frac{84000 c - 38797}{6000}
=
-\frac{451}{102000}
<0
. 
\end{align*}
On $[\, 741/1603, 799/12^3)$ 
we have 
\begin{align*}
\frac{d}{dx}
\sigma^2(x)
&\le
\frac{d}{dx}
\Bigl(
x(1-x)
+
\frac{2}{5\cdot 12}(5x-2)(6-12x)
+
\frac{2}{5^2 12^2}(5^2 x-11)(67-12^2x)
\\&\qquad
+\frac{2}{5^3 12^3}(5^3x-57)(799-12^3x)
\Bigr)
+\frac{1}{2\cdot 5^3}
\\
&
=
-\frac{1512000 x - 696743}{108000}
\le 
-\frac{1512000 (741/1603) - 696743}{108000}
=
-\frac{501853}{24732000}
<0
. 
\end{align*}

On $[\, 799/12^3, 800/12^3)$, we have
$\langle 12^3 x\rangle=12^3x-799$, 
$\langle 5^3 x\rangle=5^3x-57$, 
$\langle 12^3 x\rangle-\langle 5^3 x\rangle=1603x-742$, 
$742/1603 \in [\, 799/12^3, 800/12^3)$, and
$$
V\bigl(\langle 5^3 x\rangle, \langle 12^3 x\rangle\bigr)
=
\begin{cases}
(12^3x-799)(58-5^3x) & x\in [\, 799/12^3, 742/1603),
\\
(5^3x-57)(800-12^3x) & x\in [\,  742/1603, 800/12^3).
\end{cases}
$$
On $[\, 799/12^3, 742/1603)$ 
we have 
\begin{align*}
\frac{d}{dx}
\sigma^2(x)
&\le
\frac{d}{dx}
\Bigl(
x(1-x)
+
\frac{2}{5\cdot 12}(5x-2)(6-12x)
+
\frac{2}{5^2 12^2}(5^2 x-11)(67-12^2x)
\\&\qquad\qquad
+
\frac{2}{5^3 12^3}(12^3x-799)(58-5^3x)\Bigr)
+\frac{1}{2\cdot 5^3}
\\
&
=
-\frac{1512000 x - 698471}{108000}
\le 
-\frac{1512000 (799/12^3) - 698471}{108000}
=
-\frac{109}{18000}
<0
. 
\end{align*}
On $[\,  742/1603, 800/12^3)$ 
we have the estimate
\begin{align*}
\frac{d}{dx}
\sigma^2(x)
&\le
\frac{d}{dx}
\Bigl(
x(1-x)
+
\frac{2}{5\cdot 12}(5x-2)(6-12x)
+
\frac{2}{5^2 12^2}(5^2 x-11)(67-12^2x)
\\&\qquad\qquad
+
\frac{2}{5^3 12^3}(5^3x-57)(800-12^3x)
\Bigr)
+\frac{1}{2\cdot 5^3}
\\
&
=
-\frac{378000 x - 174217}{27000}
\le 
-\frac{378000 (742/1603) - 174217}{27000}
=
-\frac{172307}{6183000}
<0
.
\end{align*}

On $[\, 800/12^3, 801/12^3)$, we have
$\langle 12^3 x\rangle=12^3x-800$, 
$\langle 5^3 x\rangle =5^3x-57$, 
$\langle 12^3 x\rangle-\langle 5^3 x\rangle =1603x-743$, 
$743/1603\in [\, 800/12^3, 801/12^3)$, and
$$
V\bigl(\langle 5^3 x\rangle, \langle 12^3 x\rangle\bigr)
=
\begin{cases}
(12^3x-800)(58-5^3x) & x\in [\, 800/12^3, 743/1603),
\\
(5^3x-57)(801-12^3x) & x\in [\,  743/1603, 801/12^3).
\end{cases}
$$
On $[\, 800/12^3, 743/1603)$ 
we have 
\begin{align*}
\frac{d}{dx}
\sigma^2(x)
&\le
\frac{d}{dx}
\Bigl(
x(1-x)
+
\frac{2}{5\cdot 12}(5x-2)(6-12x)
+
\frac{2}{5^2 12^2}(5^2 x-11)(67-12^2x)
\\&\qquad\qquad
+
\frac{2}{5^3 12^3}(12^3x-800)(58-5^3x)
\Bigr)
+\frac{1}{2\cdot 5^3}
\\
&
=
-\frac{378000 x - 174649}{27000}
\le 
-\frac{378000 (800/12^3) - 174649}{27000}
=
-\frac{13}{1000}
<0
. 
\end{align*}
On $[\,  743/1603, 801/12^3)$ 
we have 
\begin{align*}
\frac{d}{dx}
\sigma^2(x)
&\le
\frac{d}{dx}
\Bigl(
x(1-x)
+
\frac{2}{5\cdot 12}(5x-2)(6-12x)
+
\frac{2}{5^2 12^2}(5^2 x-11)(67-12^2x)
\\&\qquad\qquad
+
\frac{2}{5^3 12^3}(5^3x-57)(801-12^3x)
\Bigr)
+\frac{1}{2\cdot 5^3}
\\
&
=
-\frac{504000 x - 232331}{36000}
\le 
-\frac{504000 (743/1603) - 232331}{36000}
=
-\frac{292201}{8244000}
<0
. 
\end{align*}

\section{Type II case, $\Sigma_{17/8}$}

 Put $c = \dfrac{101}{225}$. 
By $17^{20} = 1 \mod 225$ and $8^{20} = 1 \mod 225$, 
we have
$\langle 17^{20k+n}c\rangle = \langle 17^{n}c\rangle$, 
$\langle 8^{20k+n}c\rangle = \langle 8^{n}c\rangle$, 
and
\begin{align*}
\sigma^2(c)
&=
V(c,c) 
+
2\frac{17^{20}8^{20}}{17^{20}8^{20}-1}
\sum_{n=1}^{20}
\frac1{17^n8^n}
V\bigl(\langle 17^{n}c\rangle, \langle 8^{n}c\rangle\bigr)
\\
&=
\frac{11889992279972830720767520926085417248816985132}{47443115357778029674816911148866079222136821875}
.
\end{align*}
We divide $[\,0,1/2)$ into 
 $[\,0,\afrac7{17})$, 
 $[\,\afrac7{17}, \afrac{129}{17^2})$, 
 $[\,\afrac{129}{17^2},  \afrac{2205}{17^3}\,)$, 
 $[\,\afrac{2205}{17^3}, c)$, 
 $[\,c, \afrac{130}{17^2})$, 
 $[\,\afrac{130}{17^2}, \afrac{29}{8^2})$, 
 $[\,\afrac{29}{8^2}, \afrac8{17})$, 
 $[\,\afrac{8}{17}, \afrac12)$, 
and prove $\sigma^2(x) < \sigma^2(c) $ ($x\not=c$) on each. 

\subsection{ $[\,0,\afrac7{17})$ part.}
By applying (\ref{Eq:est}) for $N=0$  we have
\begin{align*}
\sigma^2(x)-\sigma^2(c)
&\le
\Bigl(x(1-x)+ \frac1{2(8\cdot 17-1)}-\sigma^2(c)\Bigr)_{x=7/17}
\\&
=
-\frac{128815977821313849714281109828547509123006480671
}{27422120676795701152044174644044593790395083043750}
<0
.
\end{align*}

\subsection{ $[\,\afrac{8}{17}, \afrac12)$ part.}
We have $\langle 8x\rangle=8x-3$, $\langle 17x \rangle= 17x-8$, and
$\langle 17x \rangle- \langle 8x\rangle= 9x-5< 0  $.
By applying (\ref{Eq:est}) for $N=1$  we have
\begin{align*}
\sigma^2(x)-\sigma^2(c)
&\le
\Bigl(x(1-x)
+\frac2{8\cdot 17}(17x-8)(4-8x)
+\frac1{2\cdot 8\cdot 17(8\cdot 17-1)}-\sigma^2(c)\Bigr)_{x=25/51}
\\&
=
-\frac{65851143629024194326305018601436462133936467243
}{219376965414365609216353397152356750323160664350000}
<0,
\end{align*}
since the quadratic function has the axis at $x=25/51$.

\subsection{ $[\,\afrac7{17}, \afrac{129}{17^2})$ and  $[\,\afrac{29}{8^2}, \afrac8{17})$ parts.}
On $[\, 7/17, 8/17)$, 
we have $\langle 8x\rangle=8x-3$, $\langle 17x \rangle= 17x-7$, 
$\langle 17x \rangle- \langle 8x\rangle= 9x-4$, $4/9\in [\, 7/17, 8/17)$, and
\begin{equation}\label{17over8second}
V\bigl(\langle 8 x\rangle, \langle 17 x\rangle\bigr)
=
\begin{cases}
(17x-7)(4-8x) & x\in [\, 7/17, 4/9),
\\
(8x-3)(8-17x) & x\in [\, 4/9, 8/17).
\end{cases}
\end{equation}
On $[\, 7/17, 4/9)$,  by applying (\ref{Eq:est}) for $N=1$  
we have
\begin{align*}
\sigma^2(x) -\sigma^2(c) 
&\le \Bigl(x(1-x) + \frac{2}{8\cdot 17} (17x-7)(4-8x)
+\frac{1}{2\cdot8\cdot 17 \cdot 135}\Bigr)\Bigm|_{x=4/9}
- \sigma^2(c) 
\\&=-\frac{566015656938466878466436482516492249472562779}{12904527377315624071550199832491573548421215550000}
< 0,
\end{align*}            
since the quadratic function has the axis at $x=8/17> 4/9$.

On $[\,\afrac49, \afrac{129}{17^2})$,  by applying (\ref{Eq:est}) for $N=1$  we have
\begin{align*}
\sigma^2(x) -\sigma^2(c) 
&\le \Bigl(x(1-x) + \frac{2}{8\cdot 17} (8x-3)(8-17x)
+\frac{1}{2\cdot8\cdot 17 \cdot 135}\Bigr)\Bigm|_{x=129/17^2}
- \sigma^2(c) 
\\&=-\frac{496536775682708567340343753658122607259520620727}{63399943004751661063526131777031100843393431997150000}
< 0,
\end{align*}            
since the quadratic function has the axis
at $x= \afrac{61}{136}\in (\afrac{129}{17^2},\afrac{29}{8^2})$. 
On $[\,\afrac{28}{8^2}, \afrac8{17}\,]$, it also implies 
\begin{align*}
\sigma^2(x) -\sigma^2(c) 
&\le \Bigl(x(1-x) + \frac{2}{8\cdot 17} (8x-3)(8-17x)
+\frac{1}{2\cdot8\cdot 17 \cdot 135}\Bigr)\Bigm|_{x=29/8^2}
- \sigma^2(c) 
\\&=-\frac{188828818692708659475201175773172089218806462049}
{3303559008592799762316851157117842828395831180800000}
< 0.
\end{align*}            

\subsection{ $[\,\afrac{130}{17^2}, \afrac{29}{8^2})$ part.}
We have
$\langle 8^2 x\rangle= 8^2 x -28$, 
$\langle 17 x\rangle= 17^2 x -130$, and
$\langle 17 x\rangle-\langle 8^2 x\rangle= 225x-102< 0$. 
By applying (\ref{Eq:est}) for $N=1$  we have
\begin{align*}
\sigma^2(x)- \sigma^2(c)
&
\le
\Bigl(x(1-x) + \frac{2}{8\cdot 17} (8x-3)(8-17x)
+ \frac{2}{8^2\cdot 17^2}(17^2 x -130)(29-8^2 x)
\Bigm|_{x=130/17^2}
\\&\qquad
+ \frac{1}{2\cdot 8^2\cdot 17^2 \cdot 135} - \sigma^2(c)
\\&=
-\frac{13128815469512828265962413880996134147948571587691
}{507199544038013288508209054216248806747147455977200000
}
<0, 
\end{align*}
since the quadratic function has the axis at $x=41589/92480< 130/17^2$.

\subsection{ $[\,\afrac{129}{17^2},  \afrac{2205}{17^3}\,)$  part.}
On $[\,129/17^2, 130/17^2)$, 
we have
$\langle 17^2 x\rangle = 17^2 x - 129$, 
$\langle 8^2 x\rangle = 8^2 x - 28$, 
$\langle 17^2 x\rangle -\langle 8^2 x\rangle = 225x-101$, and
\begin{equation}\label{17over8third}
V\bigl(\langle 8^2 x\rangle, \langle 17^2 x\rangle\bigr)
=
\begin{cases}
(17^2x-129)(29-8^2x) & x\in [\, 129/17^2, c),
\\
(8^2x-28)(130-17^2x) & x\in [\, c, 130/17^2).
\end{cases}
\end{equation}

On $[\,\afrac{129}{17^2}, \afrac{1975}{4401}\,]\subset[\, 129/17^2, c)$, 
by applying (\ref{Eq:est}) for $N=2$  
 we see                                       
\begin{align*}
\sigma^2(x) -\sigma^2(c) 
&\le \Bigl(x(1-x) + \frac{2}{8\cdot 17} (8x-3)(8-17x)
+ \frac{2}{8^2\cdot 17^2} (17^2x-129)(29-8^2x) 
\\&\qquad\qquad
+\frac{1}{2\cdot 8^2\cdot 17^2 \cdot 135}\Bigr)\Bigm|_{x=1975/4401}
- \sigma^2(c) 
\\& =-\frac{4490781390573907663663371846275625557701080906011}
{46629012752754238970154347271527731994688445528921200000}
< 0,
\end{align*}            
since the  quadratic function has the  axis  at 
$x=\afrac{8305}{18496} > \afrac49$. 

On $[\,\afrac{1975}{4401}, \afrac{2205}{17^3}\,)$, we have
$\langle 17^3 x\rangle = 17^3 x-2204$, $\langle 8^3 x\rangle = 8^3 x -229$, and
$\langle 17^3 x\rangle - \langle 8^3 x\rangle=4401x-1975\ge0$.
By applying (\ref{Eq:est}) for $N=3$  we have
\begin{align*}
\sigma^2(x) -\sigma^2(c)
&\le
\Bigl(
x(1-x) + \frac{2}{8\cdot 17} (8x-3)(8-17x)
+ \frac{2}{8^2\cdot 17^2} (17^2 x- 129)(29-8^2x)
\\&\qquad
+ \frac{2}{8^3\cdot 17^3} (8^3 x-229)(2205 - 17^3 x)
\Bigr)\Bigm|_{x=1975/4401}
+ \frac{1}{2\cdot 8^3\cdot 17^3 \cdot 135}-\sigma^2(c)
\\&
=
-\frac{967401418193550946022738868278421900957869917239371
}{6341545734374576499940991228927771551277628591933283200000
}
<0,
\end{align*}
since the  quadratic function has the
 axis at $x= \afrac{7901437}{17608192}< \afrac{1975}{4401}$.

\subsection{ $[\,\afrac{2205}{17^3}, c)$ and $[\,c, \afrac{130}{17^2})$ parts.} 
On $[\,2205/17^3, c\,]$, we have
$\langle 17^3 x\rangle = 17^3 x-2205$, $\langle 8^3 x\rangle = 8^3 x -229$, 
$\langle 17^3 x\rangle < \langle 8^3 x\rangle$. 
By applying (\ref{Eq:dest}) for $N=3$  
we have 
\begin{align*}
\frac{d}{dx}\sigma^2(x) \ge 
&\frac{d}{dx} 
\Bigl(
x(1-x) + \frac{2}{8\cdot 17} (8x-3)(8-17x)
+ \frac{2}{8^2\cdot 17^2} (17^2 x- 129)(29-8^2x)
\\&\qquad\qquad\qquad
+ \frac{2}{8^3\cdot 17^3} (17^3 x -2205)(230-8^3 x)
\Bigr)-\frac{2}{7\cdot 8^3}
\\&=
-\frac{123257344 x - 55339537}{8804096}
> -\frac{123257344 c - 55339537}{8804096}
=\frac{2404081}{1980921600}> 0.
\end{align*}

On $[\, c, 130/17^2)$ by applying (\ref{Eq:dest}) for $N=2$  we have
\begin{align*}
\frac{d}{dx} \sigma^2(x)
&\le\frac{d}{dx} 
\Bigl(
x(1-x) + \frac{2}{8\cdot 17} (8x-3)(8-17x)
+ \frac{2}{8^2\cdot 17^2} (8^2x-28)(130-17^2x)
\Bigr)+ \frac{2}{7\cdot 8^2}
\\&=
\le -\frac{647360 x - 289389}{64736}
=-\frac{647360 c - 289389}{64736}
-\frac{54167}{2913120}<0.
\end{align*}

\section{Type III case, $\Sigma_{19/10}$}

Put $c = \dfrac{2879}{19^3-10^3}$. 
Since we have $19^{30}=1$, $10^{30}=1$ mod $19^3-10^3$, 
we have
$\langle 19^{30k+n}c\rangle = \langle 19^{n}c\rangle$, 
$\langle 10^{30k+n}c\rangle = \langle 10^{n}c\rangle$, 
and
\begin{align*}
\sigma^2(c)
&=
V(c,c) 
+
2\frac{19^{30}10^{30}}{19^{30}10^{30}-1}
\sum_{n=1}^{30}
\frac1{19^n10^n}
V\bigl(\langle 19^{n}c\rangle, \langle 10^{n}c\rangle\bigr)
\\
&=
\frac{659906978895377949815569584725014615595601472145099845291352120114337795340}
{2637143544549129191950553455511356263269204426999999999999999999999988557373}.
\end{align*}

We divide $[\,0,1/2)$ into 
$[\,0, 8/19)$, 
$[\,8/19, 9/19)$, 
$[\,9/19,176/19^2)$, 
$[\,176/19^2, 177/19^2)$, \break
$[\,177/19^2,128/261)$, 
$[\,128/261, 3369/19^3)$,  
$[\,3369/19^3, 3370/19^3)$, 
$[\,3370/19^3, c)$, 
$[\,c, 3371/19^3)$, \break
$[\,3371/19^3, 178/19^2)$, 
$[\,178/19^2,1/2)$, 
and prove $\sigma^2(x) < \sigma^2(c) $ ($x\not=c$) on each. 

\subsection{$[\,0, 8/19)$ part.}

On $[\,0, 8/19)$, by applying (\ref{Eq:est}) for $N=0$  we have
\begin{align*}
&\sigma^2(x) -\sigma^2(c) 
< 
x(1-x) +\frac{1}{2\cdot 189}-\sigma^2(c) \biggm|_{x=8/19}
\\&
=
 - 
\frac{7278491220322831728019229898606825727606279513762088300356230722553923993855}
{1904017639164471276588299594879199222080365596293999999999999999999991738423306}
\\&<0.
\end{align*}

\subsection{$[\,8/19, 9/19)$ part.}
On $[\,8/19, 9/19)$, 
we have
$\langle 19 x\rangle= 19x-8$, 
$\langle 10 x\rangle= 10x-4$, and 
$$
V\bigl(\langle 10 x\rangle, \langle 19 x\rangle\bigr)
=
\begin{cases}
(19x-8)(5-10x) & x\in [\,8/19, 4/9),
\\
(10x-4)(9-19x) & x\in [\,4/9, 9/19).
\end{cases}
$$
On $[\,8/19, 4/9)$, by applying (\ref{Eq:est}) for $N=1$  we have
\begin{align*}
&\sigma^2(x) -\sigma^2(c)\le
x(1-x) +\frac{2}{10\cdot 19}(19x-8)(5-10x)
+\frac{1}{2\cdot 189\cdot 19\cdot 10}\biggm|_{4/9}-\sigma^2(c)
\\
&=
-\frac{710422023286243951366873538262085174754948112137941210713805643449447220303}
{1002114546928669092941210313094315380042297682259999999999999999999995651801740}
\\
&<0, 
\end{align*}
since the quadratic function  has the axis at $x=9/19> 4/9$. 

On $[\,4/9, 9/19)$, by applying (\ref{Eq:est}) for $N=1$  
 we have
\begin{align*}
&\sigma^2(x) -\sigma^2(c)\le
x(1-x) +\frac{2}{10\cdot 19}(10x-4)(9-19x)
+\frac{1}{2\cdot 189\cdot 19\cdot 10}\biggm|_{87/190}-\sigma^2(c)
\\
&=
-\frac{7910806567828266531744219401421129814389495662052207508905768063848855059757}
{47600440979111781914707489871979980552009139907349999999999999999999793460582650}
\\
&<0, 
\end{align*}
since  the quadratic function has the  axis at $x=87/190$. 

\subsection{$[\,9/19,176/19^2)$ and $[\,178/19^2,1/2)$ parts.}

On $[\,9/19, 1/2)$ we have
$\langle 19 x\rangle= 19x-9$, 
$\langle 10 x\rangle= 10x-4$, 
$\langle 19 x\rangle-\langle 10 x\rangle = 9x-5< 0$, 
and 
\begin{equation}\label{Eq;19/10second}
V\bigl(\langle 10 x\rangle, \langle 19 x\rangle\bigr) =
(19x-9)(5-10x)
\quad x\in [\,9/19, 1/2).
\end{equation}
On $[\,9/19,176/19^2)$, by applying (\ref{Eq:est}) for $N=1$  we have 
\begin{align*}
&\sigma^2(x) -\sigma^2(c)\le
x(1-x) +\frac{2}{10\cdot 19}(19x-9)(5-10x)
+\frac{1}{2\cdot 189\cdot 19\cdot 10}\biggm|_{176/19^2}-\sigma^2(c)
\\
&=
-\frac{217737941675476992486001240062929639295660101494138764285992908419778682413937}
{6873503677383741308483761537513909191710119802621339999999999999999970175708134660}
\\
&<0 ,
\end{align*}
since   the quadratic function has the axis at $28/57\in [\,176/19^2, 178/19^2)$. 
On $[\,178/19^2,1/2)$, it also implies
\begin{align*}
&\sigma^2(x) -\sigma^2(c)\le
x(1-x) +\frac{2}{10\cdot 19}(19x-9)(5-10x)
+\frac{1}{2\cdot 189\cdot 19\cdot 10}\biggm|_{178/19^2}-\sigma^2(c)
\\
&=
-\frac{6766458111546657129956963622021138234123747334138764285992908419779597824097}
{6873503677383741308483761537513909191710119802621339999999999999999970175708134660}
\\
&<0.
\end{align*}

\subsection{ $[\,176/19^2, 177/19^2)$ part.}

On $[\,176/19^2, 49/10^2)$, 
we have
$\langle 19^2 x\rangle = 19^2x -176$, 
$\langle 10^2 x\rangle= 10^2x-48$, and
$\langle 19^2 x\rangle-\langle 10^2 x\rangle= 261x-128< 0$.
By applying (\ref{Eq:est}) for $N=2$ we have
\begin{align*}
&\sigma^2(x)-\sigma^2(c)  
\\&
\le
x(1-x)+\frac{2}{10\cdot 19}(19x-9)(5-10x)
+\frac{2}{10^2 19^2}(19^2x -176)(49-10^2x)
\\&\qquad\qquad
+\frac{1}{2\cdot 189\cdot 19^2 10^2}
\biggm|_{x=49/10^2}
-\sigma^2(c)  
\\&
=
-\tfrac{86900923806297054512880071009286961604260222301441501781153612769777499920909}
{9520088195822356382941497974395996110401827981469999999999999999999958692116530000}
<0, 
\end{align*}
since the quadratic function has the axis at
$x=88489/180500> 49/10^2$.

On $[\,49/10^2,177/19^2)$, we have
$\langle 19^2 x\rangle = 19^2x -176$, 
$\langle 10^2 x\rangle= 10^2x-49$, and
$\langle 19^2 x\rangle-\langle 10^2 x\rangle= 261x-127> 0$. 
By applying (\ref{Eq:est}) for $N=2$  we have
\begin{align*}
&\sigma^2(x)-\sigma^2(c)  
\\&
\le
x(1-x) +\frac{2}{10\cdot 19}(19x-9)(5-10x)
+\frac{2}{10^2 19^2}(10^2x -49)(177-19^2x)
\\&\qquad\qquad
+\frac{1}{2\cdot 189\cdot 19^2 10^2}
\biggm|_{x=177/19^2}
-\sigma^2(c)  
\\&
=
-\tfrac{492244691787876371203608742557539740727579386088387642859929084197794135978023}
{68735036773837413084837615375139091917101198026213399999999999999999701757081346600}
<0,
\end{align*}
since the quadratic function has the axis at
$x=88589/180500> 177/19^2$.

\subsection{ $[\,177/19^2,128/261)$, $[\,128/261, 3369/19^3)$,  $[\,3371/19^3, 178/19^2)$ parts.}

On $[\,177/19^2, 178/19^2)$, 
we have
$\langle 19^2 x\rangle= 19^2x-177$, 
$\langle 10^2 x\rangle= 10^2x-49$, 
$\langle 19^2 x\rangle-\langle 10^2 x\rangle = 261x-128$, and 
\begin{equation}\label{Eq;19/10third}
V\bigl(\langle 10^2 x\rangle, \langle 19^2 x\rangle\bigr)
=
\begin{cases}
(19^2x-177)(50-10^2x) & x\in [\,177/19^2, 128/261),
\\
(10^2x-49)(178-19^2x) & x\in [\,128/261, 178/19^2).
\end{cases}
\end{equation}

On $[\,177/19^2, 128/261)$, by applying (\ref{Eq:est}) for $N=2$  we have
\begin{align*}
&\sigma^2(x)-\sigma^2(c)  
\\&
\le
x(1-x) +\frac{2}{10\cdot 19}(19x-9)(5-10x)
+\frac{2}{10^2 19^2}(19^2x-177)(50-10^2x)
\\&\qquad\qquad
+\frac{1}{2\cdot 189\cdot 19^2 10^2}
\biggm|_{x=128/261}
-\sigma^2(c)  
\\&
=
-\tfrac{691599624111739717411339488132671655847488610428626059959003766787660890038263}
{160127883453732034361075995929340654576958746648325399999999999999999305201400034600}
<0, 
\end{align*}
since the quadratic function has the axis at
$x=1779/3610> 128/261$. 

On $[\,128/261, 3369/19^3)$, by applying (\ref{Eq:est}) for $N=2$  
we have
\begin{align*}
&\sigma^2(x)-\sigma^2(c)  
\\&
\le
x(1-x) +\frac{2}{10\cdot 19}(19x-9)(5-10x)
+\frac{2}{10^2 19^2}(10^2x-49)(178-19^2x)
\\&\qquad\qquad
+\frac{1}{2\cdot 189\cdot 19^2 10^2}
\biggm|_{x=3369/19^3}
-\sigma^2(c)  
\\&
=
-\tfrac{3469534034151502550713590364547560800986360294871939072434399395404439079546939}
{24813348275355306123626379150425212182073532487463037399999999999999892334306366122600}
<0,
\end{align*}
since the quadratic function has the axis at
$x=88689/180500\in [\,3369/19^3, 3371/19^3)$. 
On $[\,3371/19^3, 178/19^2)$, it also implies
\begin{align*}
&\sigma^2(x)-\sigma^2(c)  
\\&
\le
x(1-x) +\frac{2}{10\cdot 19}(19x-9)(5-10x)
+\frac{2}{10^2 19^2}(10^2x-49)(178-19^2x)
\\&\qquad\qquad
+\frac{1}{2\cdot 189\cdot 19^2 10^2}
\biggm|_{x=3371/19^3}
-\sigma^2(c)  
\\&
=
-\tfrac{309938706743947299794717489787058688265275894403187814486879879080889481955879}
{4962669655071061224725275830085042436414706497492607479999999999999978466861273224520}
<0. 
\end{align*}

\subsection{ $[\,3369/19^3, 3370/19^3)\subset[\,491/10^3, 492/10^3)$ part.}
We have
$\langle 19^3 x\rangle= 19^3x-3369$, 
$\langle 10^3 x\rangle= 10^3x-491$, 
$\langle 19^3 x\rangle-\langle 10^3 x\rangle = 5859x-2878$, and 
$$
V\bigl(\langle 10^3 x\rangle, \langle 19^3 x\rangle\bigr)
=
\begin{cases}
(19^3x-3369)(492-10^3x) & x\in [\,3369/19^3, 2878/5859),
\\
(10^3x-491)(3370-19^3x) & x\in [\,2878/5859, 3370/19^3). 
\end{cases}
$$

On $[\,3369/19^3, 2878/5859)$, by applying (\ref{Eq:est}) for $N=3$  we have
\begin{align*}
&\sigma^2(x)-\sigma^2(c)  
\\&
\le
x(1-x) +\frac{2}{10\cdot 19}(19x-9)(5-10x)
+\frac{2}{10^2 19^2}(10^2x-49)(178-19^2x)
\\&\qquad
+\frac{2}{10^3 19^3}(19^3x-3369)(492-10^3x) 
\\&\qquad\qquad
+\frac{1}{2\cdot 189\cdot 19^3 10^3}
\biggm|_{x=2878/5859}
-\sigma^2(c)  
\\&
=
-\tfrac{4204933565133208234438019770636914745496299059344373435050395191822580968287}
{36176335144124954255177692302704785219526946329585999999999999999999843030042814000}
<0,
\end{align*}
since the quadratic function has the axis at $x=11797269/24006500> 2878/5859$. 

On $[\,2878/5859, 3370/19^3)$, by applying (\ref{Eq:est}) for $N=3$  we have
\begin{align*}
&\sigma^2(x)-\sigma^2(c)  
\\&
\le
x(1-x) +\frac{2}{10\cdot 19}(19x-9)(5-10x)
+\frac{2}{10^2 19^2}(10^2x-49)(178-19^2x)
\\&\qquad
+\frac{2}{10^3 19^3}(10^3x-491)(3370-19^3x)
+\frac{1}{2\cdot 189\cdot 19^3 10^3}
\biggm|_{x=23588679/48013000}
-\sigma^2(c)  
\\&
=
-\tfrac{7749017324660324949080945224318123073601627327116694362171996977022237046512004601}
{124066741376776530618131895752126060910367662437315186999999999999999461671531830613000000}
<0,
\end{align*}
since the quadratic function has the axis at $x=23588679/48013000$. 

\subsection{ $[\,3370/19^3, c)$ and $[\,c, 3371/19^3)$ parts.}
We have
$\langle 19^3 x\rangle= 19^3x-3370$, 
$\langle 10^3 x\rangle= 10^3x-491$, 
$\langle 19^3 x\rangle-\langle 10^3 x\rangle = 5859x-2879$, and 
$$
V\bigl(\langle 10^3 x\rangle, \langle 19^3 x\rangle\bigr)
=
\begin{cases}
(19^3x-3370)(492-10^3x) 
& x\in [\,3370/19^3, c),
\\
(10^3x-491)(3371-19^3x)
& x\in [\,c, 3371/19^3).
\end{cases}
$$
On $[\,3370/19^3, c)$, by applying (\ref{Eq:dest}) for $N=3$  
we have
\begin{align*}
\frac{d}{dx}
\sigma^2(x)
&\ge
\frac{d}{dx}
\Bigl(
x(1-x)
+
\frac{2}{10\cdot 19}
(19x-9)(5-10x)
+
\frac{2}{10^2 19^2}
(10^2x-49)(178-19^2x)
\\&\qquad\qquad
+
\frac{2}{10^3 19^3}
(19^3x-3370)(492-10^3x) 
\Bigr)
-\frac{2}{9\cdot 10^3}
\\
&
=
-\frac{432117000 x - 212352983}{30865500}
\ge
-\frac{432117000 c - 212352983}{30865500}
=\frac{1766419}{2870491500}
> 0.
\end{align*}
On $[\,c, 3371/19^3 )$, by applying (\ref{Eq:dest}) for $N=3$  
\begin{align*}
\frac{d}{dx}
\sigma^2(x)
&\le
\frac{d}{dx}
\Bigl(
x(1-x)
+
\frac{2}{10\cdot 19}
(19x-9)(5-10x)
+
\frac{2}{10^2 19^2}
(10^2x-49)(178-19^2x)
\\&\qquad\qquad
+
\frac{2}{10^3 19^3}
(10^3x-491)(3371-19^3x)
\Bigr)
+\frac{2}{9\cdot 10^3}
\\
&
=
-\frac{43211700 x - 21231397}{3086550}
\le 
-\frac{43211700 c - 21231397}{3086550}
=-\frac{186179}{287049150}
< 0.
\end{align*}

\section{Type V case, $\Sigma_{8/5}$}

Put $c = \dfrac{13690}{8^5-5^5}= \frac{13690}{29643}$. 
By 
$$
8^{40}= 5^{40}=1 \quad\hbox{mod}\quad 8^5-5^5
,$$
we have $V\bigl(\langle 5^{n+40k}c\rangle, \langle 8^{n+40k}c\rangle\bigr)
=V\bigl(\langle 5^{n}c\rangle, \langle 8^{n}c\rangle\bigr)$
and
\begin{align*}
\sigma^2(c)
&=
V\bigl(\langle c\rangle, \langle c\rangle\bigr)
+
2
\frac{5^{40}8^{40}}{5^{40}8^{40}-1}
\sum_{n=1}^{40}
\frac1{5^n8^n}
V\bigl(\langle 5^{n}c\rangle, \langle 8^{n}c\rangle\bigr)
\\&=
\frac{2693647024766931825274236270928683791388436386344146949339630859940359610}
{10622921229838049651870392375050239999999999999999999999999999999121292551}
.
\end{align*}

We divide $[\,0,1/2)$ into 
 $[\,0,2/5)$, 
 $[\,2/5, 29/8^2)$, 
 $[\,29/8^2,236/8^3)$,
 $[\,236/8^3,18/39)$, 
 $[\,18/39,1891/8^4)$,
 $[\,1891/8^4,15133/8^5)$,
 $[\,15133/8^5, c)$,
 $[\,c, 15134/8^5)$,
 $[\,15134/8^5, 1892/8^4)$,
 $[\,1892/8^4, 1893/8^4)$,\break
 $[\,1893/8^4, 179/387)$,
 $[\,179/387,237/8^3)$,
 $[\,237/8^3,30/8^2)$,
 $[\,30/8^2,12/5^2)$,
 $[\,12/5^2, 1/2)$, 
and prove $\sigma^2(x) < \sigma^2(c) $ ($x\not=c$) on each.

\subsection{ $[\,0,2/5)$ part.}
On $[\,0,2/5)$, by applying (\ref{Eq:est}) for $N=0$  we have
\begin{align*}
&\sigma^2(x) - \sigma^2(c)
\le 
x(1-x) + \frac{1}{2\cdot 39} - \sigma^2(c)\biggl|_{x=2/5}
\\&
=
-\frac{397731589368168741350186856696540338652588547976578236212312238894974663}
{531146061491902482593519618752511999999999999999999999999999999956064627550}
<0.
\end{align*}            

\subsection{ $[\,2/5, 29/8^2)$ and  $[\,12/5^2, 1/2)$ parts.}
On $[\,2/5,1/2)$, we have
$\langle 8x\rangle = 8x-3$, $\langle 5x\rangle = 5x-2$, 
$\langle 8x\rangle -\langle 5x\rangle = 3x-1> 0$, and 
\begin{equation}\label{Eq:8/5second}
V\bigl(\langle 5 x\rangle, \langle 8 x\rangle\bigr)= (5x-2)(4-8x)
\quad x\in [\,2/5,1/2). 
\end{equation}

On $[\,2/5,29/8^2)$, by applying (\ref{Eq:est}) for $N=1$   we have
\begin{align*}
&\sigma^2(x) - \sigma^2(c)
\le 
x(1-x) + \frac{2}{5\cdot 8}(5x-2)(4-8x) + \frac{1}{2\cdot 39\cdot 5\cdot 8}
 - \sigma^2(c)\biggl|_{x=29/8^2}
\\&
=
-\frac{101297915872141228476188260727611530712100269251206445552563093056478125169}
{217557426787083256870305635841028915199999999999999999999999999982004071444480}
<0,
\end{align*}            
since the quadratic function has the axis at $x=7/15\in  (29/8^2,12/5^2)$.
On $[\,12/5^2,1/2)$ it also implies
\begin{align*}
&\sigma^2(x) - \sigma^2(c)
\le 
x(1-x) + \frac{2}{5\cdot 8}(5x-2)(4-8x) + \frac{1}{2\cdot 39\cdot 5\cdot 8}
 - \sigma^2(c)\biggl|_{x=12/5^2}
\\&
=
-\frac{47677554184119599314426194214157347730517709595315647242462447781631054947}
{106229212298380496518703923750502399999999999999999999999999999991212925510000}
<0.
\end{align*}            

\subsection{ $[\,30/8^2,12/5^2)$ part.}
We have
$\langle 8^2x\rangle = 8^2x-30$, $\langle 5^2x\rangle = 5^2x-11$, and
$\langle 8^2x\rangle -\langle 5^2x\rangle = 39x-19< 0$ by $19/39> 12/5^2$. 
By applying (\ref{Eq:est}) for $N=2$  we have
\begin{align*}
&\sigma^2(x) - \sigma^2(c)
\le 
\biggl(
x(1-x) + \frac{2}{5\cdot 8}(5x-2)(4-8x) 
+ \frac{2}{5^2\cdot 8^2}(8^2x-30)(12-5^2x) 
\\
&\hskip10em
+ \frac{1}{2\cdot 39\cdot 5^2\cdot 8^2}
 - \sigma^2(c)\biggr)\biggl|_{x=1879/4000}
\\&
=
-\frac{8022607981398559967692644941121137833765667070501007117587983290720337355809}
{33993347935481758885985255600160767999999999999999999999999999997188136163200000}
<0,
\end{align*}            
since the quadratic function  has the axis at $x=1879/4000$.

\subsection{ $[\,29/8^2,236/8^3)$ and  $[\,237/8^3,30/8^2)$ parts.}

On $[\,29/8^2, 30/8^2)$, 
we have
$\langle 8^2x\rangle = 8^2x-29$, $\langle 5^2x\rangle = 5^2x-11$, 
$\langle 8^2x\rangle -\langle 5^2x\rangle = 39x-18$. 
Hence 
\begin{equation}\label{Eq:8/5third}
V\bigl(\langle 5^2 x\rangle, \langle 8^2 x\rangle\bigr)= 
\begin{cases}
(8x^2-29)(12-5^2x) & x\in [\,29/8^2, 18/39),
\\
(5x^2-11)(30-8^2x) & x\in [\,18/39,30/8^2).
\end{cases}
\end{equation}

On $[\,29/8^2, 236/8^3)\subset[\,29/8^2, 18/39)$, by applying (\ref{Eq:est}) for $N=2$  we have
\begin{align*}
&\sigma^2(x) - \sigma^2(c)
\le 
\biggl(
x(1-x) + \frac{2}{5\cdot 8}(5x-2)(4-8x) 
+ \frac{2}{5^2\cdot 8^2}(8^2x-29)(12-5^2x) 
\\
&\hskip10em
+ \frac{1}{2\cdot 39\cdot 5^2\cdot 8^2}
 - \sigma^2(c)\biggr)\biggl|_{x=236/8^3}
\\&
=
-\frac{124455417301813681838964979418237654242005385024128911051261861286851136751}
{4351148535741665137406112716820578303999999999999999999999999999640081428889600}
<0,
\end{align*}            
since the quadratic function  has the axis at $x=3733/8000> 236/8^3$.

On $[\,237/8^3,30/8^2)$, by applying (\ref{Eq:est}) for $N=2$  we have
\begin{align*}
&\sigma^2(x) - \sigma^2(c)
\le 
\biggl(
x(1-x) + \frac{2}{5\cdot 8}(5x-2)(4-8x) 
+ \frac{2}{5^2\cdot 8^2}(5^2x-11)(30-8^2x) 
\\
&\hskip10em
+ \frac{1}{2\cdot 39\cdot 5^2\cdot 8^2}
 - \sigma^2(c)\biggr)\biggl|_{x=237/8^3}
\\&
=
-\frac{34953429305342483790720514911642285574417232077212515364037956147975432359}
{13923675314373328439699560693825850572799999999999999999999999998848260572446720}
<0,
\end{align*}            
since the quadratic function has the axis at $x=1847/4000<  237/8^3$.

\subsection{ $[\,236/8^3,18/39)$ part.}
On $[\,236/8^3, 237/8^3)$, we have
$\langle 8^3x\rangle = 8^3x-236$, $\langle 5^2x\rangle = 5^3x-57$,
$\langle 8^3x\rangle -\langle 5^3x\rangle = 387x-179< 0$ for $x< 179/387$
and $\langle 8^3x\rangle -\langle 5^3x\rangle \ge 0$ otherwise. Hence
\begin{equation}\label{Eq:8/5fourth}
V(\langle 5^3x\rangle ,\langle 8^3x\rangle)
=
\begin{cases}
(8^3x-236)(58-5^3x) & x\in [\,236/8^3, 179/387),
\\
(5^3x-57)(237-8^3x) & x\in [\,179/387, 237/8^3).
\end{cases}
\end{equation}
On $[\,236/8^3,18/39)\subset [\,236/8^3, 179/387)$, by applying (\ref{Eq:est}) for $N=3$  we have
\begin{align*}
&\sigma^2(x) - \sigma^2(c)
\le 
\biggl(
x(1-x) + \frac{2}{5\cdot 8}(5x-2)(4-8x) 
+ \frac{2}{5^2\cdot 8^2}(8^2x-29)(12-5^2x) 
\\
&\hskip10em
+ \frac{2}{5^3\cdot 8^3}(8^3x-236)(58-5^3x) 
+ \frac{1}{2\cdot 39\cdot 5^3\cdot 8^3}
 - \sigma^2(c)\biggr)\biggl|_{x=18/39}
\\&
=
-\frac{3828368178678017519096077150055078687013990716720124844534877688229710583}
{3535308185290102924142466582416719871999999999999999999999999999707566160972800}
<0,
\end{align*}            
since the quadratic function  has the axis at $x=7447/16000> 18/39$.

\subsection{ $[\,179/387,237/8^3)$ part.}
Note $18/39< 179/387$. By applying (\ref{Eq:est}) for $N=3$  we have
\begin{align*}
&\sigma^2(x) - \sigma^2(c)
\le 
\biggl(
x(1-x) + \frac{2}{5\cdot 8}(5x-2)(4-8x) 
+ \frac{2}{5^2\cdot 8^2}(5^2x-11)(30-8^2x) 
\\
&\hskip10em
+ \frac{2}{5^3\cdot 8^3}(5^3x-57)(237-8^3x) 
+ \frac{1}{2\cdot 39\cdot 5^3\cdot 8^3}
 - \sigma^2(c)\biggr)\biggl|_{x=179/387}
\\&
=
-\frac{3313125491469560382482426194227489368358575253886025183417679871513529049889}
{1508488807984938532324481702512734240767999999999999999999999999875220730378163200}
<0,
\end{align*}            
since the quadratic function  has the axis at $x=206569/448000<179/387$.

\subsection{ $[\,18/39,1891/8^4)$ and  $[\,1893/8^4, 179/387)$ parts.}
On $[\,18/39,1891/8^4)$, by applying (\ref{Eq:est}) for $N=3$  we have
\begin{align*}
&\sigma^2(x) - \sigma^2(c)
\le 
\biggl(
x(1-x) + \frac{2}{5\cdot 8}(5x-2)(4-8x) 
+ \frac{2}{5^2\cdot 8^2}(5^2x-11)(30-8^2x) 
\\
&\hskip10em
+ \frac{2}{5^3\cdot 8^3}(8^3x-236)(58-5^3x) 
+ \frac{1}{2\cdot 39\cdot 5^3\cdot 8^3}
 - \sigma^2(c)\biggr)\biggl|_{x=1891/8^4}
\\&
=
-\frac{9726402459982296129168487419534163159067571323540024582460729840582190470101}
{22277880502997325503519297110121360916479999999999999999999999998157216915914752000}
<0, 
\end{align*}            
since the quadratic function  has the  axis at $x=51739/112000\in [\,1892/8^4, 1893/8^4)$. 
On $[\,1893/8^4, 179/387)$ it also implies
\begin{align*}
&\sigma^2(x) - \sigma^2(c)
\le 
\biggl(
x(1-x) + \frac{2}{5\cdot 8}(5x-2)(4-8x) 
+ \frac{2}{5^2\cdot 8^2}(5^2x-11)(30-8^2x) 
\\
&\hskip10em
+ \frac{2}{5^3\cdot 8^3}(8^3x-236)(58-5^3x) 
+ \frac{1}{2\cdot 39\cdot 5^3\cdot 8^3}
 - \sigma^2(c)\biggr)\biggl|_{x=1893/8^4}
\\&
=
-\frac{3437633091918170735261215133504421079067571323540024582460729841102385279909}
{22277880502997325503519297110121360916479999999999999999999999998157216915914752000}
<0.
\end{align*}            

\subsection{ $[\,1892/8^4, 1893/8^4)$ part.}
We have
$\langle 8^4 x\rangle = 8^4 x-1892$, 
$\langle 5^4 x\rangle = 5^4 x-288$, 
$\langle 8^4 x\rangle -\langle 5^4 x\rangle 
=
3471x-1604
$ and $1604/3471$ belongs to this interval. 

On $[\,1892/8^4, 1604/3471)$, by applying (\ref{Eq:est}) for $N=4$  we have
\begin{align*}
&\sigma^2(x) - \sigma^2(c)
\le 
\biggl(
x(1-x) + \frac{2}{5\cdot 8}(5x-2)(4-8x) 
+ \frac{2}{5^2\cdot 8^2}(5^2x-11)(30-8^2x) 
\\
&\hskip10em
+ \frac{2}{5^3\cdot 8^3}(8^3x-236)(58-5^3x) 
+ \frac{2}{5^4\cdot 8^4}(8^4x-1892)(289-5^4x) 
\\
&\hskip10em
+ \frac{1}{2\cdot 39\cdot 5^4\cdot 8^4}
 - \sigma^2(c)\biggr)\biggl|_{x=2661121/5760000}
\\&
=
-\frac{29564673245218695476520359370042348245385024128911051261861297143402050473201}
{4351148535741665137406112716820578303999999999999999999999999999640081428889600000000}
<0,
\end{align*}            
since the quadratic function  has the axis at $x=2661121/5760000$.

On $[\,1604/3471,1893/8^4)$, by applying (\ref{Eq:est}) for $N=4$  we have
\begin{align*}
&\sigma^2(x) - \sigma^2(c)
\le 
\biggl(
x(1-x) + \frac{2}{5\cdot 8}(5x-2)(4-8x) 
+ \frac{2}{5^2\cdot 8^2}(5^2x-11)(30-8^2x) 
\\
&\hskip10em
+ \frac{2}{5^3\cdot 8^3}(8^3x-236)(58-5^3x) 
+ \frac{2}{5^4\cdot 8^4}(5^4x-288)(1893-8^4x) 
\\
&\hskip10em
+ \frac{1}{2\cdot 39\cdot 5^4\cdot 8^4}
 - \sigma^2(c)\biggr)\biggl|_{x=1604/3471}
\\&
=
-\frac{698779331598069751150440137934495732570641016504944855635230157214305752517907}
{5600635227136581052426495559864567621222399999999999999999999999536726312213109760000}
<0,
\end{align*}            
since the quadratic function  has the axis at $x=10641013/23040000< 1604/3471$.

\subsection{ $[\,1891/8^4,15133/8^5)$ and  $[\,15134/8^5, 1892/8^4)$ parts.}
On $[\,1891/8^4,1892/8^4)$, we have
$\langle 8^4 x\rangle = 8^4 x-1891$, 
$\langle 5^4 x\rangle = 5^4 x-288$, 
$\langle 8^4 x\rangle -\langle 5^4 x\rangle 
=
3471x-1603
$,  $1603/3471\in(15133/8^5, c)$, 
and
\begin{equation}\label{Eq:8/5fifth}
V(\langle 5^4x\rangle ,\langle 8^4x\rangle)
=
\begin{cases}
(8^4x-1891)(289-5^4x) & x\in [\,1891/8^4, 1603/3471),
\\
(5^4x-288)(1892-8^4x) & x\in [\,1603/3471, 1892/8^4).
\end{cases}
\end{equation}

On $[\,1891/8^4,15133/8^5)$, by applying (\ref{Eq:est}) for $N=4$  we have
\begin{align*}
&\sigma^2(x) - \sigma^2(c)
\le 
\biggl(
x(1-x) + \frac{2}{5\cdot 8}(5x-2)(4-8x) 
+ \frac{2}{5^2\cdot 8^2}(5^2x-11)(30-8^2x) 
\\
&\hskip10em
+ \frac{2}{5^3\cdot 8^3}(8^3x-236)(58-5^3x) 
+ \frac{2}{5^4\cdot 8^4}(8^4x-1891)(289-5^4x)
\\
&\hskip10em
+ \frac{1}{2\cdot 39\cdot 5^4\cdot 8^4}
 - \sigma^2(c)\biggr)\biggl|_{x=15133/8^5}
\\&
=
-\frac{63426067042528236105120343586887175381622823532807866387433549238510198566647}
{7128921760959144161126175075238835493273599999999999999999999999410309413092720640000}
<0,
\end{align*}            
since the quadratic function  has the axis at $x= 1182651/2560000> 15133/8^5$.

On $[\,15134/8^5,1892/8^4)$, by applying (\ref{Eq:est}) for $N=4$  we have
\begin{align*}
&\sigma^2(x) - \sigma^2(c)
\le 
\biggl(
x(1-x) + \frac{2}{5\cdot 8}(5x-2)(4-8x) 
+ \frac{2}{5^2\cdot 8^2}(5^2x-11)(30-8^2x) 
\\
&\hskip10em
+ \frac{2}{5^3\cdot 8^3}(8^3x-236)(58-5^3x) 
+ \frac{2}{5^4\cdot 8^4}(5^4x-288)(1892-8^4x)
\\
&\hskip10em
+ \frac{1}{2\cdot 39\cdot 5^4\cdot 8^4}
 - \sigma^2(c)\biggr)\biggl|_{x=15134/8^5}
\\&
=
-\frac{11317873665183752312518460754481578805405705883201966596858387310002977399247}
{1782230440239786040281543768809708873318399999999999999999999999852577353273180160000}
<0
\end{align*}            
since the quadratic function has the axis at $x=886699/1920000 < 15134/8^5$.

\subsection{ $[\,15133/8^5, c)$ and $[\,c, 15134/8^5)$ parts.}
We have
$\langle 5^5x \rangle = 5^5x-1443$, 
$\langle 8^5x \rangle = 8^5x-15133$, and
\begin{align*}\label{Eq:8/5sixth}
V
\bigl(\langle 5^5 x\rangle, \langle 8^5 x\rangle\bigr)
&= 
\begin{cases}
(8^5x-15133)(1444-5^5x)& x \in [\,15133/8^5, c),
\\
(5^5x-1443)(15134-8^5x)& x \in [\,c, 15134/8^5). 
\end{cases}
\end{align*}
On $[\,1603/3471, c)$, by applying (\ref{Eq:dest}) for $N=5$  we have
\begin{align*}
\frac{d}{dx}
\sigma^2(x)
&\ge
\frac{d}{dx}
\Bigl(
x(1-x)
+
\frac{2}{5\cdot 8}(5x-2)(4-8x)
+
\frac{2}{5^2 8^2}(5^2x-11)(30-8^2x)
\\&\qquad\qquad
+
\frac{2}{5^3 8^3}(8^3x-236)(58-5^3x)
+
\frac{2}{5^4 8^4}(5^4x-288)(1892-8^4x)
+
\\&\qquad\qquad
\frac{2}{5^5 8^5}(8^5x-15133)(1444-5^5x)
\Bigr)
-\frac{2}{4\cdot 5^5}
\\
&
=
-\frac{225280000 x - 104042989}{10240000}
\ge 
-\frac{225280000 c - 104042989}{10240000}
=\frac{63122927}{303544320000}> 0.
\end{align*}

On $[\,15133/8^5, 1603/3471)$, by applying (\ref{Eq:dest}) for $N=5$  we have
\begin{align*}
\frac{d}{dx}
\sigma^2(x)
&\ge
\frac{d}{dx}
\Bigl(
x(1-x)
+
\frac{2}{5\cdot 8}(5x-2)(4-8x)
+
\frac{2}{5^2 8^2}(5^2x-11)(30-8^2x)
\\&\qquad\qquad
+
\frac{2}{5^3 8^3}(8^3x-236)(58-5^3x)
+
\frac{2}{5^4 8^4}(8^4x-1891)(289-5^4x)
\\&\qquad\qquad
+
\frac{2}{5^5 8^5}(8^5x-15133)(1444-5^5x)
\Bigr)
-\frac{2}{4\cdot 5^5}
\\
&
=
-\frac{225280000 x -  104070757}{10240000}
\ge 
-\frac{225280000 \cdot\frac{1603}{3471} -  104070757}{10240000}
=\frac{105757547}{35543040000}> 0.
\end{align*}

On $[\,c,15134/8^5)$, by applying (\ref{Eq:dest}) for $N=5$  we have
\begin{align*}
\frac{d}{dx}
\sigma^2(x)
&\le
\frac{d}{dx}
\Bigl(
x(1-x)
+
\frac{2}{5\cdot 8}(5x-2)(4-8x)
+
\frac{2}{5^2 8^2}(5^2x-11)(30-8^2x)
\\&\qquad\qquad
+
\frac{2}{5^3 8^3}(8^3x-236)(58-5^3x)
+
\frac{2}{5^4 8^4}(5^4x-288)(1892-8^4x)
+
\\&\qquad\qquad
\frac{2}{5^5 8^5}(5^5x-1443)(15134-8^5x)
\Bigr)
+\frac{2}{4\cdot 5^5}
\\
&
=
-\frac{563200000 x - 260100843}{25600000}
\le 
-\frac{563200000 c - 260100843}{25600000}
=-\frac{38710951}{758860800000}< 0.
\end{align*}

%%%%%%%%%%%%%%%%%%%%%%%%%%%%%%%%%%%%%%%%%
%%%%%% the references %%%%%%%%%%%%%%%%%%
%%%%%%%%%%%%%%%%%%%%%%%%%%%%%%%%%%%%%%%%%

\end{document}